\documentclass[a4paper,reqno]{amsart}

\usepackage{amssymb}
\usepackage{amstext}
\usepackage{amsmath}
\usepackage{amscd}
\usepackage{amsthm}
\usepackage{amsfonts}
\usepackage{enumerate}
\usepackage{graphicx}
\usepackage{latexsym}
\usepackage[all]{xy}
\input xy
\xyoption{all}

\usepackage{pstricks,comment}

\newtheorem{theorem}{Theorem}[section]

\newtheorem{corollary}[theorem]{Corollary}
\newtheorem{lemma}[theorem]{Lemma}
\newtheorem{proposition}[theorem]{Proposition}
\newtheorem{definition-theorem}[theorem]{Definition-Theorem}
\newtheorem{question}[theorem]{Question}

\theoremstyle{definition}
\newtheorem{definition}[theorem]{Definition}
\newtheorem{remark}[theorem]{Remark}
\newtheorem{example}[theorem]{Example}

\newcommand{\add}{\operatorname{add}\nolimits}

\newcommand{\pd}{\operatorname{proj.dim}\nolimits}

\newcommand{\bo}{\operatorname{b}\nolimits}
\newcommand{\CM}{\operatorname{CM}\nolimits}

\newcommand{\per}{\operatorname{per}\nolimits}
\newcommand{\Sl}{\operatorname{SL}\nolimits}
\newcommand{\Gl}{\operatorname{GL}\nolimits}
\newcommand{\Image}{\operatorname{Im}\nolimits}
\newcommand{\Ker}{\operatorname{Ker}\nolimits}

\newcommand{\End}{\operatorname{End}\nolimits}
\newcommand{\Ext}{\operatorname{Ext}\nolimits}
\newcommand{\Hom}{\operatorname{Hom}\nolimits}
\newcommand{\HHom}{\mathcal{H}\strut\kern-.2em\operatorname{om}\nolimits}
\newcommand{\C}{\mathsf{C}}
\newcommand{\D}{\mathsf{D}}
\newcommand{\K}{\mathsf{K}}
\newcommand{\AAA}{{\mathcal A}}
\newcommand{\BB}{{\mathcal B}}
\newcommand{\CC}{{\mathcal C}}
\newcommand{\PP}{{\mathcal P}}
\newcommand{\SSS}{{\mathcal S}}
\newcommand{\TT}{{\mathcal T}}

\def\add{\operatorname{add}}

\def\proj{\operatorname{proj}}
\def\sgn{\operatorname{sgn}\nolimits}

\def\Z{\Bbb Z}
\def\m{\mathfrak m}
\def\n{\mathfrak n}
\def\p{\mathfrak p}
\def\P{\mathfrak P}
\def\Tr{\operatorname{Tr}}
\def\mod{\operatorname{mod}\nolimits}
\def\modz{\operatorname{mod^\Z}}
\def\flz{\operatorname{fl^\Z}}
\def\lmodz{\operatorname{\underline{mod}^{\Z}}}
\def\cmz{\operatorname{CM^\Z}}
\def\lcmz{\operatorname{\underline{CM}^{\Z}}}

\def\latz{\operatorname{CM^\Z_0}}
\def\llatz{\operatorname{\underline{CM}^\Z_0}}

\def\Hom{\operatorname{Hom}}
\def\lhom{\operatorname{\underline{Hom}}}
\def\lend{\operatorname{\underline{End}}}

\def\RHom{\operatorname{{\bf R}Hom}}
\def\Tor{\operatorname{Tor}}

\def\m{\mathfrak m}
\def\n{\mathfrak n}
\def\p{\mathfrak p}

\def\k{k}
\def\l{\ell}

\begin{document}
\title{Tilting and cluster tilting for quotient singularities}
\author{Osamu Iyama and Ryo Takahashi}
\address{O. Iyama: Graduate School of Mathematics, Nagoya University, Chikusa-ku, Nagoya, 464-8602 Japan}
\email{iyama@math.nagoya-u.ac.jp}
\thanks{The first author was supported by JSPS Grant-in-Aid for Scientific Research 21740010 and 21340003.}
\address{R. Takahashi: Department of Mathematical Sciences
Faculty of Science, Shinshu University,
3-1-1 Asahi, Matsumoto, Nagano 390-8621, Japan}
\email{takahasi@math.shinshu-u.ac.jp}
\thanks{The second author was supported by JSPS Grant-in-Aid for Scientific Research 22740008.}
\dedicatory{Dedicated to Professor Shiro Goto on the occasion of his 65th birthday}

\begin{abstract}
We shall show that the stable categories of graded Cohen-Macaulay
modules over quotient singularities have tilting objects.
In particular, these categories are triangle equivalent to derived categories of finite dimensional algebras.
Our method is based on higher dimensional Auslander-Reiten theory, which gives cluster tilting objects in the stable categories of (ungraded) Cohen-Macaulay modules.
\end{abstract}
\maketitle
\tableofcontents

The aim of this paper is to discuss tilting theoretic aspects of representation theory of Cohen-Macaulay rings.
Tilting theory is a generalization of Morita theory, and it has been fundamental in representation theory of associative algebras.
While Morita theory realizes abelian categories as module categories over rings,
tilting theory realizes triangulated categories as derived categories over rings.
The key role is played by tilting objects in triangulated categories (Definition \ref{tilting}), and the tilting theorem due to Happel, Rickard and Keller asserts that
an algebraic triangulated category having a tilting object is triangle equivalent to a derived category of a ring (Theorem \ref{equivalence}).

In 1970s Auslander and Reiten initiated the representation theory of Cohen-Macaulay modules over orders, which generalize both finite dimensional algebras and commutative rings.
The theory is based on the fundamental notions of almost split sequences, Auslander-Reiten duality and Auslander algebras,
and has been developed by a number of commutative and non-commutative algebraists (see the books \cite{ARS,ASS,CR1,CR2,Y}).
One of the finest situations is given by Kleinian singularities $R$: Auslander and Herzog showed that $R$ is representation-finite
(i.e. $R$ has only finitely many isomorphism classes of indecomposable Cohen-Macaulay modules), and that
the Auslander-Reiten quiver of $R$ coincides with the McKay quiver of the corresponding group (algebraic McKay correspondence).

Recently not only representation theorists \cite{Ar,LP} but also people working on Kontsevich's homological mirror symmetry conjecture \cite{KST1,KST2,T,U} have studied
the stable categories of graded Cohen-Macaulay modules over Gorenstein singularities, especially hypersurface singularities of dimension two.
They proved that for certain singularities the stable categories contain tilting objects,
and in particular that they are triangle equivalent to derived categories of some finite dimensional algebras.

Now it is natural to ask the following:

\begin{question}
Let $R$ be a graded Gorenstein ring. When does the stable category of graded Cohen-Macaulay $R$-modules contain a tilting object?
\end{question}

The aim of this paper is to show the existence of tilting objects for quotient singularities of arbitrary dimension (Theorem \ref{tiltobj}).
Our method is to apply cluster tilting theory to show that a certain naturally constructed object $U$ is tilting.
Cluster tilting theory is one of the most active areas in recent representation theory which is closely related to the notion of Fomin-Zelevinsky cluster algebras.
It has an aspect of higher dimensional analogue of Auslander-Reiten theory, which is based on the notion of higher almost split sequences and higher Auslander algebras.
It has introduced a quite new perspective in the representation theory of Cohen-Macaulay modules \cite{BIKR,DH,KR1,KR2,KMV,I1,I2,IR,IY,IW,TV}.
A typical example is given by quotient singularities of arbitrary dimension $d$.
They have $(d-1)$-cluster tilting subcategories of the stable categories of Cohen-Macaulay modules (Theorem \ref{cluster}).
In particular this means in the case $d=2$ that Kleinian singularities are representation-finite.
For quotient singularities, the skew group algebras are higher Auslander algebras,
and the Koszul complexes are higher almost split sequences.
They will play a crucial role in the proof of our main theorem.

The following is a picture we have in mind:
\[\begin{array}{|c||c|c|c|}
\hline
&\mbox{generator}&\mbox{finite dimensional algebra}&\mbox{commutative ring}\\ \hline\hline
\mbox{Tilting theory}&\mbox{tilting object}&\mbox{derived category}&\underline{\CM}^{\Z}(R)\\ \hline
\mbox{Cluster tilting theory}&\mbox{cluster tilting object}&\mbox{cluster category}&\underline{\CM}(R)\\ \hline
\end{array}\]
This suggests studying the connection between cluster categories and stable categories of (ungraded) Cohen-Macaulay modules.
This will be done in \cite{AIR} for cyclic quotient singularities.

The authors make this paper as self-contained as possible for the reader who is not an expert.
In particular, we give detailed proofs for many important `folklore' results, such as graded Auslander-Reiten duality (Theorem \ref{310}), the description of $\End_R(S)$
as a skew group ring (Theorem \ref{endomorphism ring})
and the tilting theorem for algebraic triangulated categories (Theorem \ref{equivalence}).

We refer to \cite{ARS,ASS} for general background materials in non-commutative algebras, and to \cite{BH,Ma,Y} for ones in commutative algebras.

\medskip\noindent
{\bf Conventions }
All modules in this paper are left modules. The composition $fg$ of morphisms (respectively, arrows) means first $g$ next $f$.
For a Noetherian ring $R$, we denote by $\mod(R)$ the category of finitely generated $R$-modules, and by $\proj(R)$ the category of finitely generated projective $R$-modules.
If moreover $R$ is a $\Z$-graded ring, we denote by $\modz(R)$ the category of finitely generated $\Z$-graded $R$-modules.

For an additive category $\AAA$ and a set ${\bf M}$ of objects in $\AAA$, we denote by $\add {\bf M}$ the full subcategory of $\AAA$ consisting of modules which are isomorphic to direct summands of finite direct sums of objects in ${\bf M}$.
For a single object $M\in\AAA$, we simply write $\add M$ instead of $\add\{ M\}$.
When $\AAA$ is the category $\mod(R)$ (respectively, $\modz(R)$), we often denote $\add{\bf M}$ by $\add_R{\bf M}$ (respectively, $\add^{\Z}_R{\bf M}$).

\medskip\noindent
{\bf Acknowledgements }
The authors express their gratitude to R.-O. Buchweitz, who suggested constructing tilting objects by using syzygies. This led them to our main Theorem \ref{tiltobj}.
They are also grateful to Y. Yoshino and H. Krause for valuable suggestions on skew group algebras and the tilting theorem for algebraic triangulated categories.
Results in this paper were presented in Nagoya (June 2008), Banff (September 2008), Sherbrooke (October 2008), Lincoln (November 2008), Kyoto (November 2008) and Osaka (November 2009).
The authors thank the organizers of these meetings and seminars.

\section{Our results}

\subsection{Preliminaries}
Let us recall the notions of tilting objects \cite{R} and silting objects \cite{KV,AI} which play a crucial role in this paper.

\begin{definition}\label{tilting}
Let $\TT$ be a triangulated category.
\begin{itemize}
\item[(a)] For an object $U\in\TT$, we denote by $\mathsf{thick}(U)$ the smallest full triangulated subcategory of $\TT$ containing $U$ and
closed under isomorphism and direct summands.
\item[(b)] We say that $U\in\TT$ is \emph{tilting} (respectively, \emph{silting}) if 
$\Hom_{\TT}(U,U[n])=0$ for any $n\neq0$ (respectively, $n>0$) and $\mathsf{thick}(U)=\TT$.
\end{itemize}
\end{definition}

For example, for any ring $\Lambda$, the homotopy category
$\K^{\bo}(\proj\Lambda)$ of bounded complexes of finitely generated
projective $\Lambda$-modules has a tilting object $\Lambda$.
Moreover a certain converse of this statement holds for a triangulated category which is \emph{algebraic} (see Appendix for the definition).
Namely we have the following tilting theorem in algebraic triangulated categories \cite{Ke1} (cf. \cite{Bo}),
where an additive functor $F:\TT\to\TT'$ of additive categories is called an \emph{equivalence up to direct summands}
if it is fully faithful and any object $X\in\TT'$ is isomorphic to a direct summand of $FY$ for some $Y\in\TT$.

\begin{theorem}\label{equivalence}
Let $\TT$ be an algebraic triangulated category with a tilting object $U$.
Then there exists a triangle equivalence $\TT\to\K^{\bo}(\proj\End_{\TT}(U))$
up to direct summands.
\end{theorem}

For the convenience of the reader, we shall give an elementary proof in Appendix.

\medskip
Let $R$ be a commutative Noetherian ring which is Gorenstein and $\Z$-graded.
We assume that $(R,\m)$ is \emph{${}^\ast\!$local} in the sense that the set of graded proper ideals of $R$ has a unique maximal element $\m$.
We denote by $\modz(R)$ the category of finitely generated $\Z$-graded $R$-modules.
For $X,Y\in\modz(R)$, the morphism set $\Hom^{\Z}_R(X,Y)$ in $\modz(R)$ consists of homogeneous homomorphisms of degree $0$.
Then $\Hom^{\Z}_R(X,Y(i))$ consists of homogeneous homomorphisms of degree $i$, and we have
\[\Hom_R(X,Y)=\bigoplus_{i\in\Z}\Hom^{\Z}_R(X,Y(i)).\]
We call $X\in\modz(R)$ \emph{graded Cohen-Macaulay} if the following equivalent conditions are satisfied (e.g. \cite[(1.5.6), (2.1.17) and (3.5.11)]{BH}).
\begin{itemize}
\item $\Ext^i_R(X,R)=0$ for any $i>0$.
\item $X_\m$ is a maximal Cohen-Macaulay $R_\m$-module (i.e. ${\rm depth} X_\m=d$ or $X_\m=0$).
\item $\Ext^i_{R_\m}(X_\m,R_\m)=0$ for any $i>0$.
\item $X_\p$ is a maximal Cohen-Macaulay $R_\p$-module for any prime ideal $\p$ of $R$.
\item $\Ext^i_{R_\p}(X_\p,R_\p)=0$ for any $i>0$ and any prime ideal $\p$ of $R$.
\end{itemize}
We denote by $\cmz(R)$ the category of graded Cohen-Macaulay $R$-modules.

We denote by $\lmodz(R)$ the \emph{stable category} of $\modz(R)$ \cite{ABr}. Thus $\lmodz(R)$ has the same objects as $\modz(R)$,
and the morphism set is given by
\[\lhom^{\Z}_R(X,Y):=\Hom^{\Z}_R(X,Y)/P^{\Z}(X,Y)\]
for any $X,Y\in\lmodz(R)$, where $P^{\Z}(X,Y)$ is the submodule of $\Hom^{\Z}_R(X,Y)$ consisting of morphisms which factor through graded projective $R$-modules.
For a full subcategory $\CC$ of $\modz(R)$, we denote by $\underline{\CC}$ the corresponding full subcategory of $\lmodz(R)$.
The following fact is well-known (see Appendix).

\begin{proposition}
Let $R$ be a graded Gorenstein ring. Then $\cmz(R)$ is a Frobenius category and the stable category $\lcmz(R)$ is an algebraic triangulated category.
\end{proposition}

In particular we can apply the tilting Theorem \ref{equivalence} for $\lcmz(R)$.

Recall that an additive category is called \emph{Krull-Schmidt} if any object is isomorphic to a finite direct sum of objects whose endomorphism rings are local.
Krull-Schmidt categories are important since every object can uniquely be decomposed into indecomposable objects up to isomorphism.
Another important property of $\lcmz(R)$ is the following, which will be shown at the beginning of Section \ref{Proof of Main Results}.

\begin{proposition}\label{Krull-Schmidt}
Assume that $R_0$ is Artinian. Then the categories $\cmz(R)$ and $\lcmz(R)$ are Krull-Schmidt.
\end{proposition}

\subsection{Our results}

Throughout this paper, let $k$ be a field of characteristic zero and let $G$ be a finite subgroup of $\Sl_d(k)$.
We regard the polynomial ring $S:=k[x_1,\cdots,x_d]$ as a $\Z$-graded $k$-algebra by putting $\deg x_i=1$ for each $i$.
Since the action of $G$ on $S$ preserves the grading, the invariant subalgebra $R:=S^G$ forms a $\Z$-graded Gorenstein $k$-subalgebra of $S$.
We denote by $\m$ (respectively, $\n$) the graded maximal ideal of $R$ (respectively, $S$).

Throughout this paper, we assume that $R$ is an \emph{isolated singularity} (i.e. $R_\p$ is regular for any prime ideal $\p\neq\m$ of $R$,
or equivalently, any graded prime ideal $\p\neq\m$ of $R$).
This is equivalent to saying that $G$ acts freely on $k^d\backslash\{0\}$ \cite[(8.2)]{IY}.
One important property of the category $\lcmz(R)$ is the following graded Auslander-Reiten duality.

\begin{theorem}\label{serre}
There exists a functorial isomorphism
\[\lhom^{\Z}_R(X,Y)\simeq D\lhom^{\Z}_R(Y,X(-d)[d-1])\]
for any $X,Y\in\lcmz(R)$.
\end{theorem}

This says that the triangulated category $\lcmz(R)$ has a Serre functor $(-d)[d-1]$ in the sense of Bondal-Kapranov \cite{BK}.
The graded Auslander-Reiten duality appears in \cite{AR1} in the proof of existence theorem of almost split sequences.
For the convenience of the reader we shall give in Section 2 a complete proof for graded Gorenstein isolated singularities
(see Corollary \ref{graded AR duality} and Proposition \ref{a-invariant of invariant subring}).

Our first main result is the following.

\begin{theorem}\label{main}
The $R$-module
\[T:=\bigoplus_{p=0}^{d-1}S(p)\]
is a silting object in $\lcmz(R)$. Moreover $\lhom^{\Z}_R(T,T[n])=0$ holds for any $n<2-d$.
\end{theorem}

More strongly, we shall show that $\lhom_R^\Z(S,S(i)[n])=0$ holds if $(n,i)$ does not belong to the following shadow areas (Proposition \ref{SS}).
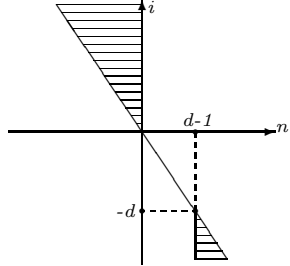
\[\begin{picture}(100,100)
\put(0,50){\vector(1,0){100}}
\put(50,0){\vector(0,1){100}}
\put(100,50){\it\scriptsize n}
\put(52,95){\it\scriptsize i}
\put(70,20){\circle*{2}}
\put(50,20){\circle*{2}}
\put(70,50){\circle*{2}}
\put(82,2){\line(-2,3){64}}
\put(50,20){\dashbox{2}(20,30)}
\put(65,52){\it\scriptsize d-1}
\put(40,18){\it\scriptsize -d}
\put(70,20){\line(0,-1){18}}
\linethickness{.1pt}
\put(70,17){\line(1,0){2}}
\put(70,14){\line(1,0){4}}
\put(70,11){\line(1,0){6}}
\put(70,8){\line(1,0){8}}
\put(70,5){\line(1,0){10}}
\put(70,2){\line(1,0){12}}
\put(50,53){\line(-1,0){2}}
\put(50,56){\line(-1,0){4}}
\put(50,59){\line(-1,0){6}}
\put(50,62){\line(-1,0){8}}
\put(50,65){\line(-1,0){10}}
\put(50,68){\line(-1,0){12}}
\put(50,71){\line(-1,0){14}}
\put(50,74){\line(-1,0){16}}
\put(50,77){\line(-1,0){18}}
\put(50,80){\line(-1,0){20}}
\put(50,83){\line(-1,0){22}}
\put(50,86){\line(-1,0){24}}
\put(50,89){\line(-1,0){26}}
\put(50,92){\line(-1,0){28}}
\put(50,95){\line(-1,0){30}}
\put(50,98){\line(-1,0){32}}
\end{picture}\]

When $d=2$, Theorem \ref{main} implies that $T=S\oplus S(1)$ is a tilting object in $\lcmz(R)$, and so $\lcmz(R)$ is triangle equivalent to $\K^{\bo}(\proj\lend^{\Z}_R(T))$.
Moreover we shall show that $\lend^{\Z}_R(T)$ is Morita equivalent to the path algebra of disjoint union of two Dynkin quivers (Example \ref{d=2}).
Thus we recover a result due to Kajiura, Saito and A. Takahashi \cite{KST1}.
Also Theorem \ref{main} gives an analogue of a result of Ueda \cite{U} based on Orlov's theorem \cite{O}, where cyclic quotient singularities with different grading are treated.


\medskip
When $d>2$, the above $T$ is not necessarily a tilting object (see Example \ref{not tilting}).
Instead of $T$, we shall give a tilting object in $\lcmz(R)$ by using syzygies as follows.

For a graded $S$-module $X$, we denote by $\Omega_SX$ the kernel of the graded projective cover.
For a graded $R$-module $X$, we denote by $[X]_{\CM}$ the maximal direct summand of $X$ which is a graded Cohen-Macaulay $R$-module.
Our second main result is the following.

\begin{theorem}\label{tiltobj}
The $R$-module
$$
U:=\bigoplus_{p=1}^d[\Omega^p_Sk(p)]_{\rm CM}
$$
is a tilting object of $\lcmz(R)$. In particular, we have a triangle equivalence
\[\lcmz(R)\to\K^{\bo}(\proj\lend^{\Z}_R(U)).\]
\end{theorem}

We shall give explicit descriptions of the algebras $\lend^{\Z}_R(T)$ and $\lend^{\Z}_R(U)$.
Let $V:=S_1$ be the degree 1 part of $S$. Let
\[E=\bigoplus_{i\ge0}E_i:=\bigoplus_{i\ge0}\bigwedge^iDV\]
be the exterior algebra of the dual vector space $DV$ of $V$.

\begin{definition}\label{definition of algebra}
\begin{itemize}
\item[(a)] We define $k$-algebras $S^{(d)}$ and $E^{(d)}$ by
\[S^{(d)}:=\left[\begin{array}{cccccc}
S_0&  0&  0&\cdots&0&0\\
S_1&S_0&  0&\cdots&0&0\\
S_2&S_1&S_0&\cdots&0&0\\
\vdots&\vdots&\vdots&\ddots&\vdots&\vdots\\
S_{d-2}&S_{d-3}&S_{d-4}&\cdots&S_0&  0\\
S_{d-1}&S_{d-2}&S_{d-3}&\cdots&S_1&S_0
\end{array}\right],\ \ \ \ \ 
E^{(d)}:=\left[\begin{array}{cccccc}
E_0&E_1&E_2&\cdots&E_{d-2}&E_{d-1}\\
0  &E_0&E_1&\cdots&E_{d-3}&E_{d-2}\\
0  &  0&E_0&\cdots&E_{d-4}&E_{d-3}\\
\vdots&\vdots&\vdots&\ddots&\vdots&\vdots\\
0&0 &0         &\cdots&E_0&E_1\\
0&0 &0         &\cdots&0  &E_0
\end{array}\right].\]
\item[(b)] Let $A$ be a $k$-algebra and $G$ a group acting on $A$ (from left).
We define the \emph{skew group algebra} $A*G$ as follows:
As a $k$-vector space, $A*G=A\otimes_kkG$. The multiplication is given by
\[(a\otimes g)(a'\otimes g')=ag(a')\otimes gg'\]
for any $a,a'\in A$ and $g,g'\in G$.
Similarly, for a $k$-algebra $A'$ and a group $G$ acting on $A'$ (from right),
we define the skew group algebra $G*A'$.
\end{itemize}
\end{definition}

The left action of $G$ on $V$ gives the right action of $G$ on $DV$.
These actions induce actions of $G$ on the $k$-algebra $S^{(d)}$ from left and the $k$-algebra $E^{(d)}$ from right.
Thus we have skew group algebras $S^{(d)}*G$ and $G*E^{(d)}$.
Notice that they are Koszul dual to each other.
Define idempotents of $kG$ by
\begin{equation}\label{idempotent}
e:=\frac{1}{\#G}\sum_{g\in G}g\ \mbox{ and }\ e':=1-e.
\end{equation}
These are idempotents of $S^{(d)}*G$ and $G*E^{(d)}$ since $kG$ is a subalgebra of $S^{(d)}*G$ and $G*E^{(d)}$ respectively.
Now we have the following descriptions of the endomorphism algebras of $T$ and $U$.

\begin{theorem}\label{End of U}
We have isomorphisms
\begin{eqnarray*}
\lend^{\Z}_R(T)&\simeq&(S^{(d)}*G)/\langle e\rangle,\\
\lend^{\Z}_R(U)&\simeq&e'(G*E^{(d)})e'
\end{eqnarray*}
of $k$-algebras, where we denote by $\langle e\rangle$ the two-sided ideal generated by $e$.
Moreover these algebras have finite global dimension.
\end{theorem}

As an immediate consequence, we have the following description of $\lcmz(R)$ as a derived category of a finite dimensional algebra.

\begin{corollary}
We have triangle equivalences:
$$
\lcmz(R)\to\K^{\bo}(\proj e'(G*E^{(d)})e')\simeq\D^{\bo}(\mod e'(G*E^{(d)})e').
$$
\end{corollary}

%

Let us observe another consequence of Theorem \ref{End of U}.

\begin{definition}
Let $\TT$ be a triangulated category.
An object $X\in\TT$ is called \emph{exceptional} if $\Hom_{\TT}(X,X[n])=0$ for any $n\neq0$ and $\End_{\TT}(X)$ is a division ring.
A sequence $(X_1,\cdots,X_m)$ of exceptional objects in $\TT$ is called an \emph{exceptional sequence} if $\Hom_{\TT}(X_i,X_j[n])=0$ for any $1\le j<i\le m$ and $n\in\Z$.
An exceptional sequence is called \emph{strong} if $\Hom_{\TT}(X_i,X_j[n])=0$ for any $1\le i,j\le m$ and $n\neq0$,
and called \emph{full} if $\mathsf{thick}(\bigoplus_{i=1}^mX_i)=\TT$.
\end{definition}

Clearly an exceptional sequence $(X_1,\cdots,X_m)$ is full and strong if and only if $\bigoplus_{i=1}^mX_i$ is a tilting object in $\TT$.
Thus the previous theorems yield the result below.

\begin{corollary}
\begin{itemize}
\item[(a)] There exists an ordering in the isomorphism classes of indecomposable direct summands of $T$ which forms a full exceptional sequence in $\lcmz(R)$.
\item[(b)] There exists an ordering in the isomorphism classes of indecomposable direct summands of $U$ which forms a full strong exceptional sequence in $\lcmz(R)$.
\end{itemize}
\end{corollary}

\section{Graded Auslander-Reiten duality}

In this section, we study the graded version of Auslander-Reiten duality.
For the basic definitions and facts on graded rings, we refer to \cite[\S 1.5 and \S 3.6]{BH}.

Throughout this section, let $R$ be a commutative Noetherian ring of Krull dimension $d$ which is Gorenstein and graded.
Moreover we assume that $(R,\m,k)$ is a ${}^\ast\!$local $k$-algebra.
We denote by $(-)^\ast$ the $R$-dual functor
\[\Hom_R(-,R):\modz(R)\to\modz(R).\]
For $X\in\modz(R)$, take a graded free resolution
\begin{equation}\label{free resolution}
\cdots\to F_1\xrightarrow{f}F_0\to X\to0.
\end{equation}
We put $\Omega X:=\Image f$. This gives the \emph{syzygy functor}
\[\Omega:\lmodz(R)\xrightarrow{}\lmodz(R).\]
Applying the functor $(-)^*$, we define $\Tr X\in\modz(R)$ by the exact sequence
\begin{equation}\label{free resolution2}
0\to X^*\to F_0^*\xrightarrow{f^*}F_1^*\to\Tr X\to0.
\end{equation}
This gives the \emph{Auslander-Bridger transpose duality}
\[\Tr:\lmodz(R)\xrightarrow{\sim}\lmodz(R).\]
Note that $\Omega^2\Tr(-)\simeq(-)^\ast$ as functors from $\lmodz(R)$ to itself.

We denote by $\flz(R)$ the category of graded $R$-modules of finite length.
Then $X\in\modz(R)$ belongs to $\flz(R)$ if and only if $X_\p=0$ for any prime ideal $\p\neq\m$
(or equivalently, any graded prime ideal $\p\neq\m$).
We denote by
\[D=\Hom_k(-,k):\flz(R)\xrightarrow{\sim}\flz(R)\]
the graded Matlis duality. 
Since $R$ is Gorenstein, there exists an integer $\alpha$ such that $\omega:=R(\alpha)$ satisfies an isomorphism
\begin{equation}\label{matlis}
\Ext^i_R(-,\omega)\simeq\left\{\begin{array}{cc}
D&(i=d),\\
0&(i\neq d)
\end{array}\right.
\end{equation}
of functors $\flz(R)\to\flz(R)$.
The integer $\alpha$ is called the \emph{$a$-invariant} of $R$ and the module $\omega$ is called the \emph{${}^\ast\!$canonical module} of $R$.
Note that $-\alpha$ is often called the \emph{Gorenstein parameter} of $R$.
In particular we have $\Ext^d_R(k,R)\simeq k(-\alpha)$.

The restriction of $\Omega$ gives an equivalence
\[\Omega:\lcmz(R)\xrightarrow{\sim}\lcmz(R),\]
which gives the suspension functor $[-1]$ of the triangulated category $\lcmz(R)$.
We define the graded version of the \emph{Auslander-Reiten translation}
\[\tau:\lcmz(R)\xrightarrow{\Omega^d\Tr}\lcmz(R)\xrightarrow{\Hom_R(-,\omega)}\lcmz(R).\]

\begin{definition}\cite{A3}
We denote by $\latz(R)$ the full subcategory of $\cmz(R)$ consisting of $N\in\cmz(R)$ such that $N_\p$ is $R_\p$-free
for any prime ideal $\p\neq\m$ (or equivalently, any graded prime ideal $\p\ne\m$).
\end{definition}

We shall show the following graded version of Auslander-Reiten duality \cite{AR1} (\cite[(3.10)]{Y}).

\begin{theorem}\label{310}
There is a functorial isomorphism
$$
D\lhom_R(M,N)\simeq\Ext_R^1(N,\tau M)
$$
of graded $R$-modules of finite length for any $M\in\modz(R)$ and $N\in\latz(R)$.
\end{theorem}

To prove this, we need the following graded version of \cite[(7.1)]{A2}\cite[(3.9)]{Y}.

\begin{lemma}\label{gprp}
There is a functorial isomorphism
\[\lhom_R(M,N)\simeq\Tor_1^R(\Tr M,N)\]
of graded $R$-modules for any $M,N\in\modz(R)$.
\end{lemma}

\begin{proof}
The assignment $\phi\otimes y\mapsto(x\mapsto\phi(x)y)$ makes a functorial homomorphism $\lambda_{M,N}:M^\ast\otimes_RN\to\Hom_R(M,N)$ of graded $R$-modules.
Clearly $\lambda_{M,N}$ is an isomorphism if $M$ is a graded free $R$-module.
Applying $-\otimes_RN$ to \eqref{free resolution2} and $\Hom_R(-,N)$ to \eqref{free resolution} and comparing them, we have a commutative diagram
\[\xymatrix{
&M^*\otimes_RN\ar[r]\ar[d]^{\lambda_{M,N}}&F_0^*\otimes_RN\ar[r]\ar[d]_{\wr}^{\lambda_{F_0,N}}&F_1^*\otimes_RN\ar[d]_{\wr}^{\lambda_{F_1,N}}\\
0\ar[r]&\Hom_R(M,N)\ar[r]&\Hom_R(F_0,N)\ar[r]&\Hom_R(F_1,N)
}\]
where the lower sequence is exact and the homology of the upper sequence is $\Tor_1^R(\Tr M,N)$.
Since the image of $\lambda_{M,N}$ is $P(M,N)$, the cokernel of $\lambda_{M,N}$ is $\lhom_R(M,N)$. 
Thus we have the desired isomorphism.
\end{proof}

Now we shall prove Theorem \ref{310} along the same lines as in the proof of \cite[(3.10)]{Y}.

First we remark the following:
Let $X,Y,Z\in\modz(R)$.
Then the assignment $\phi\mapsto(x\mapsto(y\mapsto\phi(x\otimes y)))$ makes a functorial isomorphism
$$
\Hom_R(X\otimes_RY,Z)\overset{\sim}{\to}\Hom_R(X,\Hom_R(Y,Z))
$$
of graded $R$-modules. This gives rise to an isomorphism
$$
\RHom_R(X\otimes_R^{\bf L}Y,Z)\simeq\RHom_R(X,\RHom_R(Y,Z))
$$
in the derived category $\D(\modz(R))$ of $\modz(R)$.

Since $\Ext_R^i(N,\omega)=0$ for $i>0$, we get isomorphisms
$$
\RHom_R(\Tr M\otimes_R^{\bf L}N,\omega)\simeq\RHom_R(\Tr M,\RHom_R(N,\omega))\simeq\RHom_R(\Tr M,\Hom_R(N,\omega))
$$
in $\D(\modz(R))$.
Thus there is a spectral sequence
$$
E_2^{i,j}=\Ext_R^i(\Tor^R_j(\Tr M,N),\omega)\Rightarrow\Ext_R^{i+j}(\Tr M,\Hom_R(N,\omega)).
$$
Since $N_\p$ is a free $R_\p$-module for any prime ideal $\p\neq\m$, we have $\Tor_j^R(\Tr M,N)_\p=\Tor_j^{R_\p}(\Tr M_\p,N_\p)=0$ for $j>0$.
Thus $\Tor_j^R(\Tr M,N)$ belongs to $\flz(R)$ for $j>0$.
By \eqref{matlis}, we have $E_2^{i,j}=0$ for $j>0$ and $i\neq d$.
On the other hand, we have $E_2^{i,j}=0$ for $i>d$ since $\omega$ has injective dimension $d$.
Consequently, we have an isomorphism
\begin{equation}\label{spectral}
\Ext_R^d(\Tor_1^R(\Tr M,N),\omega) \simeq\Ext_R^{d+1}(\Tr M,\Hom_R(N,\omega)).
\end{equation}
We have the desired assertion by
\begin{eqnarray*}
&&D\underline{\Hom}_R(M,N)\stackrel{\eqref{matlis}}{\simeq}
\Ext_R^d(\underline{\Hom}_R(M,N),\omega)\stackrel{\ref{gprp}}{\simeq}
\Ext_R^d(\Tor_1^R(\Tr M,N),\omega)\\
&\stackrel{\eqref{spectral}}{\simeq}&
\Ext_R^{d+1}(\Tr M,\Hom_R(N,\omega))
\simeq \Ext_R^1(\Omega^d\Tr M,\Hom_R(N,\omega))\\
&\simeq&\Ext_R^1(N,\Hom_R(\Omega^d\Tr M,\omega))=\Ext^1_R(N,\tau M).
\end{eqnarray*}
\qed

\begin{proposition}\label{tau for gorenstein}
There is an isomorphism
\[\tau\simeq\Omega^{2-d}(\alpha)\]
of functors $\lcmz(R)\to\lcmz(R)$.
\end{proposition}

\begin{proof}
We only prove the proposition in the case $d\ge2$, for it is easily proved in the cases $d=0,1$.
There is an exact sequence
$$
0 \to \Omega^{d-2}(M^\ast) \to F_{d-3} \to \cdots \to F_0 \to M^\ast \to 0
$$
of graded $R$-modules where each $F_i$ is free.
All modules appearing in this exact sequence are in $\cmz(R)$ because their localizations at $\m$ are maximal Cohen-Macaulay $R_\m$-modules.
Applying $\Hom_R(-,\omega)\simeq(-)^\ast(\alpha)$, we obtain an exact sequence
$$
0 \to M^{\ast\ast}(\alpha) \to F_0^\ast(\alpha) \to \cdots \to F_{d-3}^\ast(\alpha) \to \tau M \to 0
$$
of graded $R$-modules, and there is an isomorphism $M^{\ast\ast}\simeq M$.
This shows $\tau M\simeq\Omega^{2-d}M(\alpha)$.
\end{proof}

Our graded Auslander-Reiten duality gives the Serre duality of the triangulated category $\llatz(R)$.

\begin{corollary}\label{graded AR duality}
There exists a functorial isomorphism
$$
\lhom^{\Z}_R(M,N)\simeq D\lhom^{\Z}_R(N,M(\alpha)[d-1])
$$
for any $M,N\in\latz(R)$. In other words, the category $\llatz(R)$ has a Serre functor $(\alpha)[d-1]$.
\end{corollary}

\begin{proof}
We have only to take the degree zero parts of the graded isomorphisms
$$
D\lhom_R(M,N) \stackrel{\ref{310}}{\simeq} \Ext_R^1(N,\tau M)
\stackrel{\ref{tau for gorenstein}}{\simeq} \Ext_R^1(N,\Omega^{2-d}M(\alpha))
\simeq\lhom_R(N,M(\alpha)[d-1]).
$$
\end{proof}

%
%
%

\section{Skew group algebras}

Throughout this section, let $k$ be a field of characteristic zero and $G$ a finite subgroup of $\Gl_d(k)$.
Let $S:=k[x_1,\cdots,x_d]$ be a polynomial algebra and $R:=S^G$ be an invariant subalgebra.
The aim of this section is to give a description of $\End_R(S)$ as a skew group algebra $S*G$ (see Definition \ref{definition of algebra}).
Recall that we regard $S$ and $R$ as graded algebras by putting $\deg x_i=1$ for each $i$.
We have a morphism
\[\phi:S*G\to\End_R(S),\quad s\otimes g\mapsto (t\mapsto sg(t))\]
of graded algebras, where we regard $\End_R(S)$ and $S*G$ as graded algebras 
by putting $\End_R(S)_i:=\Hom^{\Z}_R(S,S(i))$ and $(S*G)_i:=S_i\otimes kG$.

We will show that $\phi$ is an isomorphism under the following assumption on $G$.

\begin{definition}
We say that $g\in\Gl_d(k)$ is a \emph{pseudo-reflection} if the rank of $g-1$ is at most one.
We call $G$ \emph{small} if $G$ does not contain a pseudo-reflection except the identity.
\end{definition}

For example, it is easily checked that any finite subgroup of $\Sl_d(k)$ is small.

We have the following main result in this section.

\begin{theorem}\label{endomorphism ring}
Let $G$ be a finite small subgroup of $\Gl_d(k)$.
Then the map $\phi:S*G\to\End_R(S)$ is an isomorphism of graded algebras.
\end{theorem}

This is due to Auslander \cite{AG,A1}.
Since there seem to be no convenient reference to a proof of this statement except Yoshino's book \cite{Y} for the special case $d=2$,
we shall include a complete proof for the convenience of the reader. 

Before we start to prove, we give easy consequences.

\begin{corollary}\label{negative}
Let $G$ be a finite small subgroup of $\Gl_d(k)$.
Then we have $\Hom^{\Z}_R(S,S(i))=0$ for any $i<0$ and $\End^{\Z}_R(S)\simeq kG$.
\end{corollary}

\begin{proof}
By Theorem \ref{endomorphism ring} we have
$\Hom^{\Z}_R(S,S(i))=\End_R(S)_i\simeq(S*G)_i$.
Since $(S*G)_i=0$ for any $i<0$ and $(S*G)_0=kG$, we have the assertions.
\end{proof}

\begin{corollary}\label{isoclass}
We have the following mutually quasi-inverse equivalences
\[\xymatrix@C=6em{
\mod(kG)\ar@<.3ex>^{S\otimes_{kG}-}[r]&\add^{\Z}_RS\ar@<.3ex>^{\Hom^{\Z}_R(S,-)}[l]&
\mod(kG)\ar@<.3ex>^{\Hom_{kG}(-,S)}[r]&\add^{\Z}_RS\ar@<.3ex>^{\Hom^{\Z}_R(-,S)}[l].}\]
In particular, we have bijections between the isomorphism classes of simple $kG$-modules
and the isomorphism classes of indecomposable direct summands of $S\in\mod^{\Z}(R)$.
\end{corollary}

\begin{proof}
We only consider the left diagram. We have an equivalence
\[\xymatrix@C=6em{
\proj\End_R^{\Z}(S)\ar@<.3ex>^{S\otimes_{kG}-}[r]&\add^{\Z}_RS\ar@<.3ex>^{\Hom^{\Z}_R(S,-)}[l].}\]
By Corollary \ref{negative}, we have $\End_R^{\Z}(S)\simeq kG$, which is semisimple.
Thus we have $\proj\End^{\Z}_R(S)\simeq\proj kG\simeq\mod(kG)$, and the assertions follow.
\end{proof}

The following notion plays a key role in the proof of Theorem \ref{endomorphism ring}.

\begin{definition}\cite{BR}
Let $A\subset B$ be an extension of commutative rings.
We say that $B$ is a \emph{separable} $A$-algebra if $B$ is projective as a $B\otimes_AB$-module.
\end{definition}

Let us prove the following general result (cf. \cite[(10.8)]{Y}).

\begin{proposition}\label{skew and end}
Let $B$ be an integral domain, $G$ a finite group acting on $B$ and $A:=B^G$.
If $B$ is a separable $A$-algebra, then the map $\phi:B*G\to\End_A(B)$ given by $b\otimes g\mapsto(t\mapsto bg(t))$ is surjective.
\end{proposition}

\begin{proof}
Without loss of generality, we can assume that the action of $G$ on $B$ is faithful.
The multiplication map $m:B\otimes_AB\to B$ is a surjective morphism of $B\otimes_AB$-modules.
Since $B$ is a separable $A$-algebra, there exists a morphism $\iota:B\to B\otimes_AB$ of $B\otimes_AB$-modules such that $m\iota=1_B$.

Let $E:=\End_A(B)$.
For $f,f'\in E$, we define an element $\gamma(f,f')\in E$ by the composition
\[\gamma(f,f'):B\xrightarrow{\iota}B\otimes_AB\xrightarrow{f\otimes f'}B\otimes_AB\xrightarrow{m}B.\]
Then the map $\gamma: E\times E\to E$ is bilinear.
Clearly we have $\gamma(1_E,1_E)=1_E$.

For $c\in B*G$, we simply denote $\gamma(\phi(c),f)$ and $\gamma(f,\phi(c))$ by $\gamma(c,f)$
and $\gamma(f,c)$ respectively.

(i) Let $g\in G$ and $f\in\End_A(B)$.
Since both $\iota$ and $m$ are morphisms of $B\otimes_AB$-modules, we have
\begin{equation}\label{phi 0}
\gamma(g,f)(xy)=g(x)\cdot\gamma(g,f)(y)\ \mbox{ and }\ \gamma(f,g)(xy)=\gamma(f,g)(x)\cdot g(y)
\end{equation}
for any $x,y\in B$.
In particular, we have $\gamma(g,f)(x)=g(x)\cdot\gamma(g,f)(1_B)$, which means
\begin{equation}\label{phi 1}
\gamma(g,f)=\phi(\gamma(g,f)(1_B)\otimes g).
\end{equation}
On the other hand, we have
\begin{equation}\label{phi 2}
\gamma(g,1_E)(x)\cdot y\stackrel{\eqref{phi 0}}{=}\gamma(g,1_E)(xy)=\gamma(g,1_E)(yx)
\stackrel{\eqref{phi 0}}{=}g(y)\cdot\gamma(g,1_E)(x)
\end{equation}
for any $x,y\in B$.
If $g\neq1$, then there exists $y\in B$ such that $g(y)\neq y$ since the action of $G$ on $B$ is faithful.
Since $B$ is a domain, \eqref{phi 2} implies that
\begin{equation}\label{phi 3}
\gamma(g,1_E)=0\ \mbox{ for any }\ g\in G\backslash\{1\}.
\end{equation}

(ii) Assume $\Image f\subseteq A$.
Then we have
\[(B\otimes_AB\xrightarrow{f\otimes f'}B\otimes_AB\xrightarrow{m}B)=
(B\otimes_AB\xrightarrow{f\otimes 1_E}B\otimes_AB\xrightarrow{m}B\xrightarrow{f'}B').\]
This implies
\begin{equation}\label{phi 4}
\gamma(f,f')=f'\cdot\gamma(f,1_E).
\end{equation}

(iii) We shall show $f=\phi(\sum_{g\in G}\gamma(g,f)(1_B)\otimes g)$ for any $f\in\End_A(B)$.
We have
\begin{eqnarray*}
&&\phi(\sum_{g\in G}\gamma(g,f)(1_B)*g)\stackrel{\eqref{phi 1}}{=}
\sum_{g\in G}\gamma(g,f)\stackrel{{\rm}}{=}\gamma(\sum_{g\in G}g,f)\\
&\stackrel{\eqref{phi 4}}{=}&f\cdot\gamma(\sum_{g\in G}g,1_E)
\stackrel{{\rm}}{=}\sum_{g\in G}f\cdot\gamma(g,1_E)
\stackrel{\eqref{phi 3}}{=}f\cdot\gamma(1_E,1_E)=f.
\end{eqnarray*}
\end{proof}

Another key notion in the proof is the following.


\begin{definition}
Let $(A,\m,\k)\subset(B,\n,\l)$ be an extension of commutative Noetherian local rings with $\m=\n\cap A$ such that $B$ is finitely generated as an $A$-algebra.
We say that $A\subset B$ is an {\em unramified} extension if $\n=\m B$ and $\k\subset\l$ is a separable extension of fields.
\end{definition}

This is related to separability of the $A$-algebra $B$ (see \cite[(III.3.1)]{BR} for more information).

\begin{proposition}\label{unramified}
If $A\subset B$ is unramified, then it is separable.
\end{proposition}

\begin{proof}
First, let us prove that $\Omega_{B|A}=0$.
By \cite[(25.2)]{Ma}, there is an exact sequence
$$
\n/\n^2 \overset{\alpha}{\to} \Omega_{B|A}\otimes_B\l \to \Omega_{\l|A} \to 0,
$$
where $\alpha(\overline{y})=d_{B|A}(y)\otimes\overline{1}$ for $y\in\n$.
Since $\Omega_{\l|A}=\Omega_{\l|\k}=0$, the map $\alpha$ is surjective.
As $\n=\m B$, the map $\beta:\m/\m^2\otimes_AB\to\n/\n^2$ sending $\overline{x}\otimes b$ to $\overline{xb}$ for $x\in\m$ and $b\in B$, is surjective.
We have $\alpha\beta(\overline{x}\otimes b)=x\cdot d_{B|A}(b)\otimes\overline{1}=d_{B|A}(b)\otimes\overline{x}=0$, hence $\alpha\beta=0$.
Therefore $\Omega_{B|A}\otimes_B\l=0$.
Nakayama's lemma implies $\Omega_{B|A}=0$.

Now, we can prove the assertion of the proposition.
There is an exact sequence
$$
0 \to I \to B\otimes_AB \overset{\mu}{\to} B \to 0
$$
of $B\otimes_AB$-modules, where $\mu(b\otimes b')=bb'$ for $b,b'\in B$.
We have $I/I^2=\Omega_{B|A}=0$, so $I=I^2$.
Fix a prime ideal $\p$ of $B\otimes_AB$.
If $\p$ does not contain $I$, then $B_\p=0$.
Suppose $\p$ contains $I$.
We have $I_\p=I_\p^2$.
Since $B$ is a finitely generated $A$-algebra, $B\otimes_AB$ is a Noetherian ring.
Hence we can apply Nakayama's lemma to see that $I_\p=0$, which shows $B_\p\simeq(B\otimes_RB)_\p$.
Consequently, $B$ is locally free over $B\otimes_AB$.
This means that $B$ is a projective $B\otimes_AB$-module.
\end{proof}

%

Let $L|K$ be a Galois extension of fields with the Galois group $G$.
Let $(B,\P)$ be a discrete valuation ring with the quotient field $L$.
Then $A:=B\cap K$ is a discrete valuation ring with the maximal ideal $\p:=\P\cap A$.
The \emph{inertia group} of $\P$ is defined by
\[T(\P):=\{g\in G\ |\ g(x)-x\in\P\ \mbox{ for any }x\in B\}.\]
We need the following basic result in algebraic number theory \cite[(I.7.20)]{S}.

\begin{proposition}\label{inertia}
The order of $T(\P)$ is equal to the length of the $B$-module $B/\p B$.
\end{proposition}

In particular, the extension $A\subset B$ is unramified if and only if $T(\P)=\{1\}$ holds and the field extension $A/\p\subset B/\P$ is separable.
Applying this observation for our setting, we have the following crucial result (cf. \cite[(10.7.2)]{Y}).

\begin{proposition}\label{small and unramified}
Let $k$ be a field of characteristic zero and $G$ a finite subgroup of $\Gl_d(k)$.
Let $S:=k[x_1,\cdots,x_d]$ and $R:=S^G$.
Then the following conditions are equivalent.
\begin{itemize}
\item[(a)] $G$ is small.
\item[(b)] For any height one prime ideal $\P$ of $S$, the extension $R_{\P\cap R}\subset S_\P$ is unramified.
\end{itemize}
\end{proposition}

\begin{proof} 
Let $\p:=\P\cap R$.
Since the characteristic of $k$ is zero, the field extension $R_{\p}/\p R_{\p}\subset S_\P/\P S_{\P}$ is separable.
By Proposition \ref{inertia}, we have only to show that the following conditions are equivalent for $g\in G$.
\begin{itemize}
\item[(i)] $g$ is a pseudo-reflection.
\item[(ii)] There exists a height one prime ideal $\P$ of $S$ such that $g\in T(\P S_\P)$.
\end{itemize}

(i)$\Rightarrow$(ii) By changing variables,
we can assume that $g(x_i)=x_i$ for any $1\le i<d$ and $g(x_d)=\zeta x_d$ for a root $\zeta$ of unity.
Let $\P:=(x_d)$ be a height one prime ideal.
Since $g$ acts on $S/(x_d)$ trivially, we have that $g\in T(\P S_\P)$.

(ii)$\Rightarrow$(i) Since $S$ is a unique factorization domain, we can write $\P=(z)$ for $z\in\mathfrak{n}$.
Since $g$ acts on $S/(z)$ trivially, $g$ acts on $S/((z)+\mathfrak{n}^2)$ trivially.
Let $\mathfrak{n}:=(x_1,\cdots,x_d)\subset S$ and $V:=S_1=\mathfrak{n}/\mathfrak{n}^2$.
The image of the map $g-1:V\to V$ is contained in $((z)+\mathfrak{n}^2)/\mathfrak{n}^2\subset V$, which is one dimensional over $k$.
Thus $g$ is a pseudo-reflection.
\end{proof}

\begin{lemma}\label{height one iso}
Let $A$ be a commutative Noetherian ring of dimension $d$ such that ${\rm depth} A_{\m}=d$ for any maximal ideal $\m$ of $A$.
Let $f:X\to Y$ be a homomorphism in $\mod A$ such that $f_{\P}:X_{\P}\to Y_{\P}$ is an isomorphism for any height one prime ideal $\P$ of $A$.
Then $f^*:Y^*\to X^*$ is an isomorphism.
In particular, $f$ is an isomorphism if $X$ and $Y$ are reflexive $A$-modules.
\end{lemma}

\begin{proof}
Let $0\to K\to X\xrightarrow{f}Y\to C\to 0$ be an exact sequence.
By our assumption $\dim K$ and $\dim C$ are at most $d-2$.
By \cite[(17.1 Ischebeck)]{Ma}, we have $\Ext_A^i(K,A)=0=\Ext_A^i(C,A)$ for $i=0,1$.
Applying $(-)^*$ to the above exact sequence, we have the assertion.
\end{proof}

Now we can prove Theorem \ref{endomorphism ring}.

The map $\phi:S*G\to\End_R(S)$ is a homomorphism of reflexive $S$-modules.
By Lemma \ref{height one iso}, we have only to show that
\[\phi_\P:S_\P*G\to\End_{R_{\P\cap R}}(S_\P)\]
is an isomorphism for any height one prime ideal $\P$ of $S$.

Let $A:=R_{\P\cap R}$ and $B:=S_\P$, and let $K$ and $L$ be the quotient fields of $A$ and $B$ respectively.
Since $G$ is small, we have by Proposition \ref{small and unramified} that $A\subset B$ is an unramified extension of discrete valuation rings.
By Proposition \ref{unramified}, $B$ is a separable $A$-algebra.
By Proposition \ref{skew and end}, $\phi_\P$ is surjective.
Since $\dim_K(L*G)=(\dim_KL)^2=\dim_K\End_K(L)$, we have that $\phi_\P$ is bijective.
\qed

\section{Cluster tilting subcategories for quotient singularities}

The notion of maximal orthogonal subcategories was introduced in \cite{I1} as a domain of higher dimensional analogue of Auslander-Reiten theory, and later
renamed cluster tilting subcategories in \cite{KR1} in the context of categorification of Fomin-Zelevinsky cluster algebras.
The notion of cluster tilting subcategories is an analogue of that of tilting objects.
It was shown in \cite{I1} that the categories $\underline{\CM}(R)$ for isolated quotient singularities $R$ of dimension $d$ have $(d-1)$-cluster tilting subcategories.
In this section we give a graded version of this observation, namely we show that the categories $\lcmz(R)$ have $(d-1)$-cluster tilting subcategories.
Results in this section will be used in the proofs of our main Theorems \ref{main} and \ref{tiltobj}.

Let us start with introducing the notion of $n$-cluster tilting subcategories:

\begin{definition}
Let $\TT$ be a triangulated category and let $n\ge1$.
We say that a full subcategory $\CC$ of $\TT$ is \emph{$n$-cluster tilting} (or \emph{maximal $(n-1)$-orthogonal}) if the following conditions are satisfied:
\begin{itemize}
\item[(i)] We have
\begin{eqnarray*}
\CC=\{X\in\TT\ |\ \Hom_{\TT}(\CC,X[j])=0\ (0<j<n)\}=\{X\in\TT\ |\ \Hom_{\TT}(X,\CC[j])=0\ (0<j<n)\}.
\end{eqnarray*}
\item[(ii)] $\CC$ is a functorially finite subcategory of $\TT$ in the sense that for any $X\in\TT$,
there exist morphisms $f:Y\to X$ and $g:X\to Z$ with $Y,Z\in\CC$ such that the following are surjective.
\[(f\cdot):\Hom_{\TT}(\CC,Y)\to\Hom_{\TT}(\CC,X),\ \ \ \ \ (\cdot g):\Hom_{\TT}(Z,\CC)\to\Hom_{\TT}(X,\CC).\]
\end{itemize}
\end{definition}

The following main result in this section is a graded version of \cite[(2.5)]{I1} (see also \cite{Bu}).

\begin{theorem}\label{cluster}
Let $k$ be a field of characteristic zero and $G$ a finite subgroup of $\Sl_d(k)$.
Let $S:=k[x_1,\cdots,x_d]$ and $R:=S^G$.
We regard $S$ and $R$ as graded algebras by putting $\deg x_i=1$ for each $i$.
If $R$ is an isolated singularity, then
\[\SSS:=\add^{\Z}_R\{S(i)\ |\ i\in\Z\}\]
is a $(d-1)$-cluster tilting subcategory of $\lcmz(R)$.
\end{theorem}

We include a complete proof here for the convenience of the reader though it is parallel to \cite[(2.5)]{I1}.
Let us recall some well-known results.
The first one is a consequence of our assumption that $R$ is an isolated singularity.

\begin{lemma}\label{Hom-finite}
Let $X,Y\in\cmz(R)$ and $n\in\Z$.
There exist only finitely many $i\in\Z$ such that $\lhom^{\Z}_R(X,Y(i))\neq0$ (respectively, $\Ext^n_R(X,Y(i))_0\neq0$, where $\Ext^n_R(X,Y(i))_0$ denotes the degree zero part of $\Ext^n_R(X,Y(i))$).
\end{lemma}

\begin{proof}
We have $\lhom_R(X,Y)=\bigoplus_{i\in\Z}\lhom_R^\Z(X,Y(i))$ and $\Ext^n_R(X,Y)=\bigoplus_{i\in\Z}\Ext^n_R(X,Y(i))_0$.
Since $R$ is an isolated singularity, we have $\dim_k\lhom_R(X,Y)<\infty$ and $\dim_k\Ext^i_R(X,Y)<\infty$ for any $i>0$ by \cite[(7.6)]{A3} (cf. \cite[(3.3)]{Y}).
Thus the assertion follows.
\end{proof}


The next one is a general result on graded modules.

\begin{lemma}\label{comparison of iso}\cite[(2.9)]{AR2}\cite[(15.2.2)]{Y}
Let $X$ and $Y$ be indecomposable objects in $\modz(R)$.
Then $X\simeq Y$ as $R$-modules if and only if $X\simeq Y(i)$ for some $i\in\Z$ as graded $R$-modules.
\end{lemma}

\noindent{\sc Proof of Theorem \ref{cluster} }
For any $X\in\lcmz(R)$, there are only finitely many $i\in\Z$ such that $\lhom^{\Z}_R(X,S(i))\neq0$ (respectively, $\lhom^{\Z}_R(S,X(i))\neq0$) by Lemma \ref{Hom-finite}.
One can easily verify that $\SSS$ is a functorially finite subcategory of $\lcmz(R)$.

(i) We shall show $\lhom^{\Z}_R(S,S(i)[j])=0$ for any $i\in\Z$ and $0<j<d-1$.

We have only to show $\Ext^j_R(S,S)=0$ for any $0<j<d-1$.
Fix $0<j<d-1$ and assume that we have shown $\Ext^\ell_R(S,S)=0$ for any $0<\ell<j$. Let
\[0\to\Omega^j_RS\to P_{j-1}\to\cdots\to P_0\to S\to0\]
be a projective resolution of the $R$-module $S$.
Applying $\Hom_R(-,S)$, one gets an exact sequence
\[0\to\End_R(S)\to\Hom_R(P_0,S)\to\cdots\to\Hom_R(P_{j-1},S)\to\Hom_R(\Omega^j_RS,S)\to\Ext^j_R(S,S)\to0\]
of $S$-modules since $\Ext^\ell_R(S,S)=0$ for any $0<\ell<j$.
Each $\Hom_R(P_i,S)$ is a projective $S$-module, and so is $\End_R(S)$ since there is an isomorphism $\End_R(S)\simeq S*G$ by Theorem \ref{endomorphism ring}.
Taking a projective presentation of the $R$-module $\Omega^j_RS$ and applying $\Hom_R(-,S)$, we see that $\Hom_R(\Omega^j_RS,S)$ is a second syzygy of a certain $S$-module.
Thus it holds that $\pd_S\Hom_R(\Omega^j_RS,S)\le d-2$.
These observations together with the above exact sequence imply $\pd_S\Ext^j_R(S,S)\le d-1$.
By Lemma \ref{Hom-finite}, the $S$-module $\Ext^j_R(S,S)$ has finite length.
Since any non-zero finite length $S$-module has projective dimension $d$ by the Auslander-Buchsbaum equality at $\m$, we have $\Ext^j_R(S,S)=0$.

(ii) We shall show that any $X\in\lcmz(R)$ satisfying $\lhom^{\Z}_R(X,S(i)[j])=0$
for any $i\in\Z$ and $0<j<d-1$ belongs to $\SSS$.

By Lemma \ref{comparison of iso}, we have only to check $X\in\add_RS$.
We have $\Ext^j_R(X,S)=0$ for any $0<j<d-1$. Let
\[0\to\Omega_R^{d-1}X\to P_{d-2}\to\cdots\to P_{0}\to X\to0\]
be a projective resolution of the $R$-module $X$.
Applying $\Hom_R(-,S)$, we obtain an exact sequence
\[0\to\Hom_R(X,S)\to\Hom_R(P_0,S)\to\cdots\to\Hom_R(P_{d-2},S)\to\Hom_R(\Omega^{d-1}_RS,S)\]
of $S$-modules since $\Ext^j_R(X,S)=0$ for any $0<j<d-1$.
Again each $\Hom_R(P_i,S)$ is a projective $S$-module, and $\pd_S\Hom_R(\Omega^{d-1}_RS,S)\le d-2$ holds as we have shown in (i).
Consequently $\Hom_R(X,S)$ is a projective $S$-module, and we get $\Hom_R(X,S)\in\add_RS$.

Since $eS=R$ holds for the idempotent defined in \eqref{idempotent}, one sees that $R$ is a direct summand of $S$ as an $R$-module.
Thus we have $X^*\in\add_R\Hom_R(X,S)\subset\add_RS$ and $X\simeq X^{**}\in\add_RS^*$.
Since $S^*$ is a Cohen-Macaulay $S$-module, it is a projective $S$-module.
Thus $X\in\add_RS$ holds, so the desired assertion follows.

(iii) We shall show that any $X\in\lcmz(R)$ satisfying $\lhom^{\Z}_R(S,X(i)[j])=0$
for any $i\in\Z$ and $0<j<d-1$ belongs to $\SSS$.

By using graded AR duality (Corollary \ref{graded AR duality}), we have $\lhom^{\Z}_R(X,S(i)[j])=0$
for any $i\in\Z$ and $0<j<d-1$. Thus we have $X\in\SSS$ by (ii).
\qed

\medskip
In particular, we have $\underline{\CM}^{\Z}(R)=\SSS$ in the case $d=2$.
This is a graded version of a classical result due to Herzog \cite{Her} and Auslander \cite{A4}.

%


The following property of $(d-1)$-cluster tilting subcategories is the graded version of \cite[(2.2.3)]{I1}.

\begin{proposition}\label{generator2}
For any $X\in\cmz(R)$, there exists an exact sequence
\[0\to C_{d-2}\to\cdots\to C_0\to X\to0\]
in $\modz(R)$ with $C_i\in\SSS$ such that the following sequence is exact:
\[0\to\Hom_R(S,C_{d-2})\to\cdots\to\Hom_R(S,C_0)\to\Hom_R(S,X)\to0.\]
\end{proposition}

\begin{proof}
Put $X_0:=X$. Inductively, we get an exact sequence $0\to X_{i+1}\to C_i\to X_i\to0$
in $\modz(R)$ with $C_i\in\SSS$ such that
\[0\to\Hom_R(S,X_{i+1})\to\Hom_R(S,C_i)\to\Hom_R(S,X_i)\to0\]
is exact.
Using $\Ext^j_R(S,S)=0$ for any $0<j<d-1$, we have $\Ext^j_R(S,X_i)=0$ for any $0<j\le i<d-1$ inductively.
By Theorem \ref{cluster}, we have $X_{d-2}\in\SSS$.
Putting $C_{d-2}:=X_{d-2}$, we obtain a desired sequence.
\end{proof}


\section{Koszul complexes are higher almost split sequences}\label{Higher almost split sequences}

The notion of Koszul complexes is fundamental in commutative algebra.
They give projective resolutions of simple modules for regular rings $S$.
In this section we study properties of the Koszul complex of $S$ as a complex of $R$-modules for a quotient singularity $R$.
We will show that the Koszul complex of $S$ is a direct sum of higher almost split sequences for graded Cohen-Macaulay modules.
This is an analogue of a result in \cite{I2} which assumes that $S$ is a formal power series ring.
Results in this section will be used in the proof of our main Theorems \ref{main} and \ref{tiltobj}.

Let $k$ be a field of characteristic zero and $G$ a finite subgroup of $\Gl_d(k)$.
Let $S:=k[x_1,\cdots,x_d]$ and $R:=S^G$.
We regard $S$ and $R$ as graded algebras by putting $\deg x_i=1$ for each $i$.
Let $V:=S_1$ be the degree $1$ part of $S$.
We denote by
\[\mathbf{K}=(0\to S\otimes\bigwedge^dV\xrightarrow{\delta_d}S\otimes\bigwedge^{d-1}V\xrightarrow{\delta_{d-1}}\cdots
\xrightarrow{\delta_2}S\otimes V\xrightarrow{\delta_1}S\to0)\]
the Koszul complex of $S$, where $\otimes:=\otimes_k$ and the map $\delta_i:S\otimes\bigwedge^iV\to S\otimes\bigwedge^{i-1}V$ is given by
\[s\otimes(v_1\wedge\cdots\wedge v_i)\mapsto\sum_{j=1}^i(-1)^{j-1}sv_j\otimes(v_1\wedge\cdots\check{v}_j\cdots\wedge v_i).\]
We have a free resolution
\begin{equation}\label{Koszul}
\widetilde{\mathbf{K}}=(0\to S\otimes\bigwedge^dV\xrightarrow{\delta_d}S\otimes\bigwedge^{d-1}V\xrightarrow{\delta_{d-1}}\cdots\xrightarrow{\delta_2}S\otimes V\xrightarrow{\delta_1}S\to k\to0)
\end{equation}
of the simple $S$-module $k$.
We regard each $S\otimes\bigwedge^iV$ as a module over the group algebra $RG$ by
\[(rg)(s\otimes(v_1\wedge\cdots\wedge v_i))=rg(s)\otimes(g(v_1)\wedge\cdots\wedge g(v_i)).\]
We regard $V$ as a vector space of degree $1$, and so $\bigwedge^iV$ as a vector space of degree $i$. Thus $S\otimes\bigwedge^iV$ is a direct sum of copies of $S(-i)$.
Clearly \eqref{Koszul} is an exact sequence of graded $RG$-modules.

We consider the decomposition of \eqref{Koszul} as a complex of $R$-modules.
Let $V_0=k,V_1,\cdots,V_n$ be the isomorphism classes of simple $kG$-modules.
Let $d_i$ be the dimension of $V_i$ over a division ring $\End_{kG}(V_i)$.

\begin{proposition}\label{n-almost split}
As a complex of graded $R$-modules, the koszul complex \eqref{Koszul} splits into a direct sum of $d_i$ copies of
\begin{equation}\label{summand of K}
\Hom_{kG}(V_i,\widetilde{\mathbf{K}})=(0\to(\Hom_{kG}(V_i,S\otimes\bigwedge^dV)\xrightarrow{\delta_d\cdot}\cdots\xrightarrow{\delta_1\cdot}\Hom_{kG}(V_i,S)\to\Hom_{kG}(V_i,k)\to0),
\end{equation}
where $\Hom_{kG}(V_i,k)$ if $i=0$ and $\Hom_{kG}(V_i,k)=0$ otherwise.
\end{proposition}

\begin{proof}
Since $kG\simeq\bigoplus_{i=0}^nV_i^{d_i}$ as $kG$-modules, we have $\widetilde{\mathbf{K}}\simeq\bigoplus_{i=0}^{d_i}\Hom_{kG}(V_i,\widetilde{\mathbf{K}})^{d_i}$
as complexes of graded $R$-modules.
\end{proof}

%

The following main result in this section is a graded version of \cite[(6.1)]{I2} asserting that the Koszul complex of $S$ gives a higher almost split sequence of $R$.

\begin{theorem}\label{gorenstein condition}
\begin{itemize}
\item[(a)] The Koszul complex \eqref{Koszul} induces an exact sequence
\[0\to\Hom^{\Z}_R(S(i),S\otimes\bigwedge^dV)\xrightarrow{\delta_d\cdot}\cdots
\xrightarrow{\delta_2\cdot}\Hom^{\Z}_R(S(i),S\otimes V)\xrightarrow{\delta_1\cdot}
\Hom^{\Z}_R(S(i),S)\]
for any $i\in\Z$, where the right map $(\delta_1\cdot)$ is surjective if $i\neq0$.
\item[(b)] The Koszul complex \eqref{Koszul} induces an exact sequence
\[0\to\Hom^{\Z}_R(S,S(i))\xrightarrow{\cdot\delta_1}\cdots\xrightarrow{\cdot\delta_{d-1}}
\Hom^{\Z}_R(S\otimes\bigwedge^{d-1}V,S(i))\xrightarrow{\cdot\delta_d}
\Hom^{\Z}_R(S\otimes\bigwedge^dV,S(i))\]
for any $i\in\Z$, where the right map $(\cdot\delta_d)$ is surjective if $i\neq-d$.
\end{itemize}
\end{theorem}

The first part is a consequence of the following observation.

\begin{lemma}\label{iso from koszul}
We have an isomorphism 
\[\xymatrix{
\mathbf{K}\otimes kG=\ar[d]_{\wr}&(S\otimes\bigwedge^dV\otimes kG\ar[r]^(0.7){\delta_d\otimes1}\ar[d]_{\wr}^{a_d}&
\cdots\ar[r]^{\delta_2\otimes1}
&S\otimes V\otimes kG\ar[r]^{\delta_1\otimes1}\ar[d]_{\wr}^{a_1}&S\otimes kG\ar[d]_{\wr}^{a_0})\\
\Hom_R(S,\mathbf{K})=&(\Hom_R(S,S\otimes\bigwedge^dV)\ar[r]^(0.7){\delta_d\cdot}&\cdots\ar[r]^(0.3){\delta_2\cdot}
&\Hom_R(S,S\otimes V)\ar[r]^{\delta_1\cdot}&\Hom_R(S,S))
}\]
of complexes of graded $R$-modules.
\end{lemma}

\begin{proof}
We define a map
\[a_i:S\otimes\bigwedge^iV\otimes kG\to\Hom_R(S,S\otimes\bigwedge^iV)
\ \mbox{ by }\ 
s\otimes w\otimes g\mapsto(t\mapsto tg(s)\otimes g(w))\ (g\in G).\]

(i) We shall show that $a_i$ is an isomorphism.

We can write $a_i$ as a composition
\[S\otimes\bigwedge^iV\otimes kG\xrightarrow{c_i}S\otimes kG\otimes\bigwedge^iV\xrightarrow{\phi\otimes1}\End_R(S)\otimes\bigwedge^iV\xrightarrow{\sim}
\Hom_R(S,S\otimes\bigwedge^iV),\]
where $c_i$ is defined by $c_i(s\otimes w\otimes g)=s\otimes g\otimes g(w)$.
Clearly $c_i$ is an isomorphism, and so is $\phi\otimes1$ by Theorem \ref{endomorphism ring}. Thus we have the assertion.

(ii) It is easily checked that both $(\delta_i\cdot)a_i$ and $a_{i-1}(\delta_i\otimes1)$ send $s\otimes(v_1\wedge\cdots\wedge v_i)\otimes g$ to
\[t\mapsto\sum_{j=1}^i(-1)^{j+1}tg(sv_j)\otimes(g(v_1)\wedge\cdots\check{g(v_j)}\cdots\wedge g(v_i)).\]
Thus the assertion follows.
\end{proof}

The second part is a consequence of the following observation.

\begin{lemma}\label{iso from koszul2}
We have an isomorphism 
\[\xymatrix{
\mathbf{K}\otimes kG=\ar[d]_{\wr}&(S\otimes\bigwedge^dV\otimes kG\ar[r]^(0.8){\delta_d\otimes1}\ar[d]_{\wr}^{b_d}&
\cdots\ar[r]^{\delta_2\otimes1}&
S\otimes V\otimes kG\ar[r]^{\delta_1\otimes1}\ar[d]_{\wr}^{b_1}&S\otimes kG\ar[d]_{\wr}^{b_0})\\
\Hom_R(\mathbf{K},S)=&(\Hom_R(S,S)\ar[r]^(0.6){-(\cdot\delta_1)}&\cdots\ar[r]^(0.2){(-1)^{d-1}(\cdot\delta_{d-1})}&
\Hom_R(S\otimes\bigwedge^{d-1}V,S)\ar[r]^{(-1)^d(\cdot\delta_d)}&\Hom_R(S\otimes\bigwedge^dV,S))
}\]
of complexes of graded $R$-modules.
\end{lemma}

\begin{proof}
Fix an isomorphism $\iota:\bigwedge^dV\to k$.
We have a non-degenerate bilinear form
\[\iota': (\bigwedge^iV)\otimes(\bigwedge^{d-i}V)\to k,\ \ \ w\otimes w'\mapsto \iota(w\wedge w').\]
We define a map
\[b_i:S\otimes\bigwedge^iV\otimes kG\to\Hom_R(S\otimes\bigwedge^{d-i}V,S)
\ \mbox{ by }\ 
s\otimes w\otimes g\mapsto(t\otimes w'\mapsto tg(s)\iota(g(w)\wedge w')).\]

(i) We shall show that $b_i$ is an isomorphism.

We denote by $\alpha:\bigwedge^iV\to D(\bigwedge^{d-i}V)$ the corresponding isomorphism to $\iota'$.
Since we can write $b_i$ as a composition
\[S\otimes\bigwedge^iV\otimes kG\xrightarrow{c_i}S\otimes kG\otimes\bigwedge^iV\xrightarrow{\phi\otimes\alpha}\End_R(S)\otimes D(\bigwedge^{d-i}V)\xrightarrow{\sim}
\Hom_R(S\otimes\bigwedge^{d-i}V,S),\]
we have the assertion by Theorem \ref{endomorphism ring}.

(ii) Using the relation
\[\sum_{j=1}^{d+1}(-1)^{j+1}u_j\iota(u_1\wedge\cdots\check{u_j}\cdots\wedge u_{d+1})=0\]
for any $u_1,\cdots,u_{d+1}\in V$, we have that both $(\cdot\delta_{d-i+1})b_i$ and $b_{i-1}(\delta_i\otimes1)$ send $s\otimes(v_1\wedge\cdots\wedge v_i)\otimes g$ to
\begin{eqnarray*}
t\otimes(v'_1\wedge\cdots\wedge v'_{d-i+1})&\mapsto&\sum_{j=1}^{d-i+1}(-1)^{j+1}tv'_jg(s)\iota(g(v_1)\wedge\cdots g(v_i)\wedge v'_1\cdots\check{v'_j}\cdots\wedge v'_{d-i+1})\\
&&=\sum_{j=1}^i(-1)^{i+j}tg(sv_j)\iota(g(v_1)\wedge\cdots\check{g(v_j)}\cdots g(v_i)\wedge v'_1\cdots\wedge v'_{d-i+1}).
\end{eqnarray*}
Thus the assertion follows.
\end{proof}


Now we are ready to prove Theorem \ref{gorenstein condition}.
The assertions follow from Lemmas \ref{iso from koszul} and \ref{iso from koszul2}
since the complex $\mathbf{K}\otimes kG$ is clearly exact.
\qed

\medskip
We end this section by stating the following application of Theorem \ref{gorenstein condition}.

\begin{proposition}\label{a-invariant of invariant subring}
We have $\Ext^{d}_R(k,R)\simeq k(d)$ in $\modz(R)$.
In particular the $a$-invariant of $R$ is $-d$.
\end{proposition}

\begin{proof}
This lemma follows from \cite[(6.4.9) and (6.4.11)]{BH}.
As an application of the results stated above, we give a direct proof for the convenience of the reader.
Since $R$ is Gorenstein, we have $\Ext^n_R(S,R)=0$ for any $n>0$.
Applying $\Hom^{\Z}_R(-,R(-d))$ to \eqref{Koszul}, we have an exact sequence
\[\Hom^{\Z}_R(S\otimes\bigwedge^{d-1}V,R)\to\Hom^{\Z}_R(S\otimes\bigwedge^dV,R(-d))\to\Ext^d_R(k,R(-d))_0\to0.\]
Since $R(-d)$ is a direct summand of $S\otimes\bigwedge^dV$, we have $\Ext^d_R(k,R(-d))_0\neq0$
by Theorem \ref{gorenstein condition}. This implies the assertion.
\end{proof}


\section{Proof of Main Results}\label{Proof of Main Results}

Before proving our main Theorems \ref{main} and \ref{tiltobj},
we deduce Proposition \ref{Krull-Schmidt} from the following well-known result.

\begin{proposition}\label{lifting}
Let $f:A\to B$ be a surjective homomorphism of finite dimensional algebras over a field.
For any idempotent $e$ of $B$, there exists an idempotent $e'$ of $A$ such that $e=f(e')$.
\end{proposition}

\begin{proof}
This is an easy consequence of the lifting idempotents property (e.g. \cite[(I.4.4)]{ASS}).
\end{proof}

Now we are ready to prove Proposition \ref{Krull-Schmidt}.

By our assumption $\dim_kR_i<\infty$ for each $i$, we know that $\End^{\Z}_R(X)$ and $\lend^{\Z}_R(X)$ are finite dimensional $k$-algebras for any $X\in\cmz(R)$.
If a finite dimensional $k$-algebra does not have idempotents except $0$ and $1$, then it is local.
Thus it is enough to show that idempotents split in $\cmz(R)$ (respectively, $\lcmz(R)$), i.e. for any idempotent $e\in\End^{\Z}_R(X)$ (respectively, $\lend^{\Z}_R(X)$),
there exists an isomorphism $f:X\xrightarrow{\sim}Y\oplus Z$ such that $e=f^{-1}{1\ 0\choose 0\ 0}f$.

This is clear for $\End^{\Z}_R(X)$; just put $Y:=\Image e$ and $Z:=\Ker e$.
For $\lend^{\Z}_R(X)$, apply Proposition \ref{lifting} to the surjective map $\End^{\Z}_R(X)\to\lend^{\Z}_R(X)$.
\qed

\subsection{Proof of Theorem \ref{main}}
We need the following application of results in the previous section.

\begin{lemma}\label{summand of Koszul}
As a complex of graded $R$-modules, the sequence (\ref{Koszul}) has a direct summand
\begin{equation}\label{M sequence}
0\to M_{d}\xrightarrow{g_d}M_{d-1}\xrightarrow{g_{d-1}}\cdots\xrightarrow{g_2}M_1\xrightarrow{g_1}M_0\to0
\end{equation}
which is exact and satisfies the following conditions.
\begin{itemize}
\item[(a)] $M_i\in\add S(-i)$ for any $0\le i\le d$,
\item[(b)] $R\oplus M_0=S$ and $R(-d)\oplus M_d=S(-d)$,
\item[(c)] Put $N_i:=\Image g_i$. Then $N_i=[\Omega_S^ik]_{\CM}$ for any $0<i<d$.
\item[(d)] We have an exact sequence
\[0\to\Hom^{\Z}_R(S(i),M_d)\xrightarrow{g_d\cdot}\Hom^{\Z}_R(S(i),M_{d-1})\xrightarrow{g_{d-1}\cdot}
\cdots\xrightarrow{g_2\cdot}\Hom^{\Z}_R(S(i),M_1)\xrightarrow{g_1\cdot}\Hom^{\Z}_R(S(i),M_0)\]
for any $i\in\Z$, where the right map $(g_1\cdot)$ is surjective if $i\neq0$.
\item[(e)] We have an exact sequence
\[0\to\Hom^{\Z}_R(M_0,S(i))\xrightarrow{\cdot g_1}\Hom^{\Z}_R(M_1,S(i))\xrightarrow{\cdot g_2}\cdots
\xrightarrow{\cdot g_{d-1}}\Hom^{\Z}_R(M_{d-1},S(i))\xrightarrow{\cdot g_d}\Hom^{\Z}_R(M_d,S(i))\]
for any $i\in\Z$, where the right map $(\cdot g_d)$ is surjective if $i\neq-d$.
\end{itemize}
\end{lemma}

\begin{proof}
Immediate from Proposition \ref{gorenstein condition} and Theorem \ref{n-almost split}.
\end{proof}


We have the following consequence.

\begin{proposition}\label{thick}
$\mathsf{thick}(T)=\lcmz(R)$.
\end{proposition}

\begin{proof}
By (a) and (b) in Lemma \ref{summand of Koszul}, we have $S(i)\in\mathsf{thick}(T)$ for any $i\in\Z$ by using the grade shifts of the exact sequence \eqref{M sequence}.
By Proposition \ref{generator2}, we have $\mathsf{thick}(T)=\lcmz(R)$.
\end{proof}

In the rest of this section, we simply write
\[(X,Y)_i^n:=\lhom^{\Z}_R(X,Y(i)[n]).\]
We have the following immediate consequence of Theorem \ref{serre}.

\begin{lemma}\label{duality}
$(X,Y)_i^n\simeq D(Y,X)_{-d-i}^{d-1-n}$.
\end{lemma}

Theorem \ref{main} follows immediately from the following proposition.

\begin{proposition}\label{SS}
Let $n,i\in\Z$.
\begin{itemize}
\item[(a)] $(S,S)_i^n=0$ if $n>0$ and $dn+(d-1)i>0$.
\item[(b)] $(S,S)_i^n=0$ if $n<d-1$ and $dn+(d-1)i<0$.
\end{itemize}
\end{proposition}

\begin{proof}
We use the exact sequence in Lemma \ref{summand of Koszul}.

(a) By Theorem \ref{cluster}, this is true for $0<n<d-1$.
Since $\lhom_R(S,S[n])=\bigoplus_{i\in\Z}(S,S)_i^n$ is finite dimensional, this is also true for sufficiently large $i$. 
Now we take any $n\ge d-1$.
We assume that $(S,S)_{i'}^{n'}=0$ holds
if the following (i) or (ii) is satisfied.
\begin{itemize}
\item[(i)] $n>n'>0$ and $dn'+(d-1)i'>0$,
\item[(ii)] $n=n'$ and $i'>i$.
\end{itemize}

For any $2\le j\le d$, we have a triangle
\begin{equation}\label{NM}
N_j\to M_{j-1}\to N_{j-1}\to N_j[1]
\end{equation}
in $\lcmz(R)$.
Since $n\ge n-j+2$ and $d(n-j+2)+(d-1)(i+j-1)=dn+(d-1)i+(d-j+1)>dn+(d-1)i$,
we have $(S,S)_{i+j-1}^{n-j+2}=0$ and $(M_{j-1},S)_i^{n-j+2}=0$ by our assumption of the induction.
Applying $(-,S)_i^{n-j+1}$ to \eqref{NM}, we have an exact sequence
\begin{equation}\label{induction sequence}
(M_{j-1},S)_i^{n-j+1}\to(N_j,S)_i^{n-j+1}\to(N_{j-1},S)_i^{n-j+2}\to(M_{j-1},S)_i^{n-j+2}=0
\end{equation}
for any $2\le j\le d$.
Thus we have a sequence of surjections
\begin{equation}\label{surjections}
(N_d,S)_i^{n-d+1}\to(N_{d-1},S)_i^{n-d+2}\to\cdots\to(N_2,S)_i^{n-1}\to(N_1,S)_i^n=(S,S)_i^n.
\end{equation}

Assume $n>d-1$.
Since $n\ge n-d+1>0$ and $d(n-d+1)+(d-1)(i+d)=dn+(d-1)i>0$, we have $(S,S)_{i+d}^{n-d+1}=0$ and $(N_d,S)_i^{n-d+1}=0$ by the induction hypothesis.
By \eqref{surjections}, we have $(S,S)_i^n=0$.

Assume $n=d-1$. Since we have that $i\ne-d$ and that $(M_{d-1},S)_i^0\stackrel{}{\to}(N_d,S)_i^0\to0$ is exact by Lemma \ref{summand of Koszul},
it holds that $(N_{d-1},S)_i^{n-d+2}=0$ by \eqref{induction sequence} for $j=d$.
Again by \eqref{surjections}, we have $(S,S)_i^n=0$.

(b) Since $d-1-n>0$ and $d(d-1-n)+(d-1)(-d-i)=-dn-(d-1)i>0$, the assertion follows from (a) and Lemma \ref{duality}.
\end{proof}

%
%

\subsection{Proof of Theorem \ref{tiltobj}}


Let us begin with proving several propositions.

\begin{proposition}\label{thick2}
$\mathsf{thick}(U)=\lcmz(R)$.
\end{proposition}

\begin{proof}
By Proposition \ref{thick}, we have only to show $S(i)\in\mathsf{thick}(U)$ for any $0\le i<d$. 

Fix $0\le i<d$ and assume $S(j)\in\mathsf{thick}(U)$ for any $0\le j<i$.
By Lemma \ref{summand of Koszul}, there is an exact sequence
\[0\to[\Omega^i_Sk(i)]_{\CM}\to M_{i-1}(i)\to\cdots\to M_1(i)\to M_0(i)\to0\]
with $[\Omega^i_Sk(i)]_{\CM}\in\add U$ and $M_{i-1}(i),\cdots,M_1(i)\in\add\bigoplus_{j=1}^{i-1}S(j)$.
Thus we have $S(i)=R(i)\oplus M_0(i)\in\mathsf{thick}(U)$.
The assertion is proved inductively.
\end{proof}


We use the notation in Lemma \ref{summand of Koszul}.
As a matter of convenience, set $N_0=N_{d+1}=0$.
For integers $n,i$ and $0\le p,q\le d+1$, we set
$$
(p,q)^n_i:=(N_p,N_q)_i^n=\lhom^{\Z}_R(N_p,N_q(i)[n]).
$$
Let us prepare the following result.

\begin{proposition}\label{1k}
One has $(p,q)^n_i=0$ if either of the following holds.
\begin{itemize}
\item[(a)] $(1,q)^{n+t}_{i+p-1-t}=0$ for any $0\le t\le p-1$.
\item[(b)] $(1,q)^{n-t}_{i+p+t}=0$ for any $0\le t\le d-p$.
\end{itemize}
\end{proposition}

\begin{proof}
For each integer $0\le j\le d$ there is a triangle
\[N_{j+1}\to M_j\to N_j\to N_{j+1}[1]\]
in $\lcmz(R)$.
Applying $(-,q)_i^n$ to this triangle, we get an exact sequence
\begin{equation}\label{long}
(j+1,q)_i^{n-1} \to (j,q)_i^n \to (M_j,q)_i^n \to (j+1,q)_i^n \to (j,q)_i^{n+1}
\end{equation}
of $R$-modules, and $(M_j,q)_i^n$ is in $\add_R(1,q)_{i+j}^n$.

(a) By \eqref{long} and our assumption, we have a sequence of injections
\[(p,q)_i^n\to(p-1,q)_i^{n+1}\to\cdots\to(2,q)_i^{n+p-2}\to(1,q)_i^{n+p-1}=0.\]

(b) By \eqref{long} and our assumption, we have a sequence of surjections
\[0=(d,q)_i^{n-d+p}\to(d-1,q)_i^{n-d+p+1}\to\cdots\to(p+1,q)_i^{n-1}\to(p,q)_i^n.\]
\end{proof}



A dual version of Proposition \ref{1k} holds:

\begin{proposition}\label{k1}
One has $(p,q)^n_i=0$ if either of the following holds.
\begin{itemize}
\item[(a)] $(p,1)^{n-t}_{i-q+1+t}=0$ for any $0\le t\le q-1$.
\item[(b)] $(p,1)^{n+t}_{i-q-t}=0$ for any $0\le t\le d-q$.
\end{itemize}
\end{proposition}

\begin{proof}
(a) By Lemma \ref{duality}, we have $(1,p)_{-d-i+q-1-t}^{d-1-n+t}=0$ for any $0\le t\le q-1$.
By Proposition \ref{1k}(a), we have $(q,p)_{-d-i}^{d-1-n}=0$.
By Lemma \ref{duality}, we have $(p,q)_i^n=0$.

(b) is shown dually.
\end{proof}

%
%

%

%

The next lemma will play an essential role in the proof of Theorem \ref{tiltobj}.

\begin{lemma}\label{keyvnsh}
One has $(1,q)^n_i=0$ if either $n\ge\max\{1,q-i\}$ or $n\le\min\{-1,q-i-1\}$.
\end{lemma}

\begin{proof}
Let $n\ge\max\{1,q-i\}$.
Proposition \ref{k1}(b) says that we have only to prove $(1,1)_{i-q-t}^{n+t}=0$ for $0\le t\le d-q$.
Since $n\ge q-i$ and $n+t\ge n\ge 1$, we have
\[d(n+t)+(d-1)(i-q-t)=(d-1)(n+i-q)+n+t\ge n+t>0.\]
Proposition \ref{SS}(a) implies the lemma.
A similar argument shows the assertion in the case $n\le\min\{-1,q-i-1\}$.
\end{proof}

Now we can prove Theorem \ref{tiltobj}.

By Proposition \ref{thick2} it suffices to prove $(U,U)_0^n=0$ for any $n\ne 0$.
Since $U\simeq\bigoplus_{p=1}^dN_p(p)$ in $\lcmz(R)$, there is an isomorphism $(U,U)_0^n\simeq\bigoplus_{1\le p,q\le d}(p,q)_{q-p}^n$.

Fix an integer $t\ge 0$.
When $n>0$, we have $n+t\ge\max\{1,q-(q-1-t)\}$.
It follows from Lemma \ref{keyvnsh} that $(1,q)_{q-1-t}^{n+t}=0$.
When $n<0$, we have $n-t\le\min\{-1,q-(q+t)-1\}$.
Again Lemma \ref{keyvnsh} implies that $(1,q)_{q+t}^{n-t}=0$.
Therefore, Proposition \ref{1k} shows that $(p,q)_{q-p}^n=0$ holds for every integer $n\ne 0$.
This completes the proof of the statement that $U$ is a tilting object.

Now we show the second statement of Theorem \ref{tiltobj}.
By Theorem \ref{equivalence}, we have a triangle equivalence $\lcmz(R)\to\K^{\bo}(\proj\lend^{\Z}_R(U))$ up to direct summands.
This is an equivalence since $\lcmz(R)$ is Krull-Schmidt by Proposition \ref{Krull-Schmidt}.
\qed

\section{Endomorphism algebras of $T$ and $U$}

In this section, we will give descriptions of the endomorphism algebras $\End_R^{\Z}(T)$ and $\End_R^{\Z}(U)$.
In particular, we shall prove Theorem \ref{End of U}.

Let $V_0=k,V_1,\cdots,V_n$ be the isomorphism classes of simple $kG$-modules.
Let $d_i$ be the dimension of $V_i$ over a division ring $\End_{kG}(V_i)$.
Let $e_i$ be the central idempotent of $kG$ corresponding to $V_i$.
Then we have $e_i=\sum_{j=1}^{d_i}e_{ij}$ for primitive idempotents $e_{ij}$ of $kG$.

The endomorphism algebra $\End_R^{\Z}(T)$ is rather easy to calculate.

\begin{theorem}\label{End of T}
We have the following isomorphisms of $k$-algebras:
\begin{eqnarray*}
\End^{\Z}_R(T)&\simeq&S^{(d)}*G,\\
\lend^{\Z}_R(T)&\simeq&(S^{(d)}*G)/\langle e\rangle.
\end{eqnarray*}
\end{theorem}

\begin{proof}
By Theorem \ref{endomorphism ring}, there exist isomorphisms
\[\Hom^{\Z}_R(S(p),S(q))\simeq\End_R(S)_{q-p}\simeq(S*G)_{q-p}=S_{q-p}*G.\]
Thus we have the isomorphism $\End^{\Z}_R(T)\simeq S^{(d)}*G$.

Notice that $\langle e\rangle$ consists of endomorphisms of $T$ factoring through $eT$.
Since $eT=\bigoplus_{i=0}^{d-1}R(i)$ is a maximal free direct summand of $T$,
we have $\langle e\rangle\subset P^{\Z}(T,T)$.
On the other hand, assume that for $0\le p,q\le d$, there exist non-zero maps
$S(q)\stackrel{a}{\to}R(r)\stackrel{b}{\to}S(p)$.
Since $R\in\add_R^{\Z}S$, we have $q\le r\le p$ by Corollary \ref{negative}.
Thus $ba$ factors through $eT$, and we have $\langle e\rangle\supset P^{\Z}(T,T)$.
\end{proof}

We give explicit information on the direct summands of $T$.
Let $e_p$ be the idempotent of $S^{(d)}$ whose $(p,p)$-entry is $1$ and the other entries are $0$. Then
\begin{equation}\label{e^p_{ij}}
\{e^p_{ij}:=e_p\otimes e_{ij}\ |\ 1\le p\le d,\ 0\le i\le n,\ 1\le j\le d_i\}
\end{equation}
gives a complete set of orthogonal primitive idempotents of $S^{(d)}*G$. We put
\begin{equation}\label{T_ip}
T_{pi}:=\Hom_{kG}(V_i,S)(p-1)\ \mbox{ for }\ (p,i)\in\{1,2,\cdots,d\}\times\{0,1,\cdots,n\}.
\end{equation}

\begin{proposition}\label{summand of T}
\begin{itemize}
\item[(a)] We have $e^p_{ij}T\simeq T_{pi}$ as graded $R$-modules.
\item[(b)] The isomorphism classes of indecomposable direct summands of $T\in\mod^{\Z}(R)$ (respectively, $T\in\lcmz(R)$) are
$T_{pi}$ for $(p,i)\in\{1,2,\cdots,d\}\times\{0,1,\cdots,n\}$ (respectively, $(p,i)\in\{1,2,\cdots,d\}\times\{1,\cdots,n\})$.
\end{itemize}
\end{proposition}

\begin{proof}
(a) We have $e^p_{ij}T\simeq e_{ij}S(p-1)\simeq\Hom_{kG}((kG)e_{ij},S(p-1))\simeq\Hom_{kG}(V_i,S)(p-1)$.

(b) This is immediate from Theorem \ref{End of T}.
\end{proof}

In the rest of this section we calculate the endomorphism algebra $\End_R^{\Z}(U)$.
We give a description of the endomorphism algebra of a bigger object
\[\widetilde{U}:=\bigoplus_{p=1}^d\Omega_S^pk(p)\in\modz(R),\]
which induces a description of $\lend^{\Z}_R(U)$. 

\begin{theorem}\label{End of tildeU}
We have the an isomorphism $\End^{\Z}_R(\widetilde{U})\simeq G*E^{(d)}$ of $k$-algebras.
\end{theorem}

Before proving Theorem \ref{End of tildeU}, we give explicit information on the direct summands of $U$.
Let $e_p$ be the idempotent of $E^{(d)}$ whose $(p,p)$-entry is $1$ and the other entries are $0$. Then
\[\{e^p_{ij}:=e_{ij}\otimes e_p\ |\ 1\le p\le d,\ 0\le i\le n,\ 1\le j\le d_i\}\]
gives a complete set of orthogonal primitive idempotents of $G*E^{(d)}$. We put
\begin{equation}\label{U_ip}
U_{pi}:=\Hom_{kG}(V_i,\Omega_S^pk)(p)\ \mbox{ for }\ (p,i)\in\{1,2,\cdots,d\}\times\{0,1,\cdots,n\}
\end{equation}
which is $\Image(\delta_p\cdot)(p)$ for the map $(\delta_p\cdot)$ in the complex \eqref{summand of K}.

\begin{proposition}\label{summand of U}
\begin{itemize}
\item[(a)] We have $e^p_{ij}U\simeq U_{pi}$ as graded $R$-modules.
\item[(b)] The isomorphism classes of indecomposable direct summands of $\widetilde{U}\in\mod^{\Z}(R)$ (respectively, $U\in\lcmz(R)$) are $U_{pi}$ for $(p,i)\in\{1,2,\cdots,d\}\times\{0,1,\cdots,n\}$
(respectively, $(p,i)\in\{1,2,\cdots,d\}\times\{1,\cdots,n\}$).
\end{itemize}
\end{proposition}

The idea of the proof of Theorem \ref{End of tildeU} is to calculate $\Hom^{\Z}_R(\Omega_S^qk(q),\Omega_S^pk(p))$
for any $1\le p,q\le d$ by using the Koszul complex \eqref{Koszul}
\[\mathbf{K}=(\cdots\to 0\to S\otimes\bigwedge^dV\xrightarrow{\delta_d}S\otimes\bigwedge^{d-1}V\xrightarrow{\delta_{d-1}}\cdots\xrightarrow{\delta_2}S\otimes V\xrightarrow{\delta_1}S\to0\to\cdots)\]
used in Section \ref{Higher almost split sequences}. We will construct the following isomorphisms.
\begin{equation}\label{desired isomorphisms}
kG\otimes\bigwedge^{q-p}DV\stackrel{\sim}{\to}
\Hom_{\C(\modz(R))}(\mathbf{K}(q)[-q],\mathbf{K}(p)[-p])\stackrel{\sim}{\to}
\Hom^{\Z}_R(\Omega_S^qk(q),\Omega_S^pk(p)).
\end{equation}

First we shall construct the right isomorphism of \eqref{desired isomorphisms}.

Since $\Image \delta_p=\Omega_S^pk$, we have a map
\begin{equation}\label{image of Koszul}
\Hom_{\C(\modz(R))}(\mathbf{K}(q)[-q],\mathbf{K}(p)[-p])\to\Hom^{\Z}_R(\Omega_S^qk(q),\Omega_S^pk(p)),
\end{equation}
where $\C(\modz(R))$ is the category of chain complexes over $\modz(R)$.

\begin{proposition}\label{Koszul and omega}
\begin{itemize}
\item[(1)] The map \eqref{image of Koszul} is an isomorphism.
\item[(2)] $\Hom^{\Z}_R(\Omega_S^qk(q),\Omega_S^pk(p))=0$ for any $1\le q<p\le d$.
\end{itemize}
\end{proposition}

\begin{proof}
(1) The key result in the proof is Theorem \ref{gorenstein condition}.

Fix any $f\in\Hom^{\Z}_R(\Omega_S^qk(q),\Omega_S^pk(p))$.
Using Theorem \ref{gorenstein condition} repeatedly, we have commutative diagrams
\begin{eqnarray*}
\xymatrix{
\cdots\ar[r]^(0.3){\delta_{q+2}(q)}&S\otimes\bigwedge^{q+1}V(q)\ar[r]^{\delta_{q+1}(q)}\ar[d]&S\otimes\bigwedge^qV(q)\ar[r]\ar[d]&\Omega_S^qk(q)\ar[r]\ar[d]^f&0\\
\cdots\ar[r]^(0.3){\delta_{p+2}(p)}&S\otimes\bigwedge^{p+1}V(p)\ar[r]^{\delta_{p+1}(p)}&S\otimes\bigwedge^pV(p)\ar[r]&\Omega_S^pk(p)\ar[r]&0,
}\\
\xymatrix{
0\ar[r]&\Omega_S^qk(q)\ar[r]\ar[d]^f&S\otimes\bigwedge^{q-1}V(q)\ar[r]^{\delta_{q-1}(q)}\ar[d]
&S\otimes\bigwedge^{q-2}V(q)\ar[r]^(0.7){\delta_{q-2}(q)}\ar[d]&\cdots\\
0\ar[r]&\Omega_S^pk(p)\ar[r]&S\otimes\bigwedge^{p-1}V(p)\ar[r]^{\delta_{p-1}(p)}
&S\otimes\bigwedge^{p-2}V(p)\ar[r]^(0.7){\delta_{p-2}(p)}&\cdots
}
\end{eqnarray*}
of graded $R$-modules. They give a chain homomorphism $\mathbf{K}(q)[-q]\to\mathbf{K}(p)[-p]$,
and the map \eqref{image of Koszul} is surjective.

To show that \eqref{image of Koszul} is injective, assume that the following diagram is commutative.
\[\xymatrix{
\cdots\ar[r]^(0.3){\delta_{q+2}(q)}&S\otimes\bigwedge^{q+1}V(q)\ar[r]^{\delta_{q+1}(q)}\ar[d]^{a_{q+1}}&S\otimes\bigwedge^qV(q)\ar[r]\ar[d]^{a_q}&\Omega_S^qk(q)\ar[r]\ar[d]^0&0\\
\cdots\ar[r]^(0.3){\delta_{p+2}(p)}&S\otimes\bigwedge^{p+1}V(p)\ar[r]^{\delta_{p+1}(p)}&S\otimes\bigwedge^pV(p)\ar[r]&\Omega_S^pk(p)\ar[r]&0
}\]
Since $\delta_p(p)\cdot a_q=0$, the map $a_q$ factors through $\delta_{p+1}(p)$
by Theorem \ref{gorenstein condition}.
Since $\Hom^{\Z}_R(S\otimes\bigwedge^qV(q),S\otimes\bigwedge^{p+1}V(p))=0$ by Corollary \ref{negative},
we have $a_q=0$.
Repeating similar argument, we have $a_i=0$ for any $i\ge q$.

Using the diagram below and Theorem \ref{gorenstein condition}, we have $a_i=0$ for any $i<q$ by a similar argument.
\[\xymatrix{
0\ar[r]&\Omega_S^qk(q)\ar[r]\ar[d]&S\otimes\bigwedge^{q-1}V(q)\ar[r]^{\delta_{q-1}(q)}\ar[d]^{a_{q-1}}
&S\otimes\bigwedge^{q-2}V(q)\ar[r]^(0.7){\delta_{q-2}(q)}\ar[d]^{a_{q-2}}&\cdots\\
0\ar[r]&\Omega_S^pk(p)\ar[r]&S\otimes\bigwedge^{p-1}V(p)\ar[r]^{\delta_{p-1}(p)}
&S\otimes\bigwedge^{p-2}V(p)\ar[r]^(0.7){\delta_{p-2}(p)}&\cdots
}\]
Thus the map \eqref{image of Koszul} is injective.

(2) By (1) we have only to show $\Hom_{\C(\modz(R))}(\mathbf{K}(q)[-q],\mathbf{K}(p)[-p])=0$.
Fix any chain homomorphism
\[\xymatrix{
S\otimes V(q)\ar[r]^{\delta_1(q)}\ar[d]^{a_1}&S(q)\ar[r]\ar[d]^{a_0}&0\\
S\otimes\bigwedge^{p-q+1}V(p)\ar[r]^{\delta_{p-q+1}(p)}&S\otimes\bigwedge^{p-q}V(p)
\ar[r]^{\delta_{p-q}(p)}&S\otimes\bigwedge^{p-q-1}V(p).
}\]
Since $\delta_{p-q}(p)\cdot a_0\cdot \delta_1(q)=\delta_{p-q}(p)\cdot \delta_{p-q+1}(p)\cdot a_1=0$,
we have $\delta_{p-q}(p)\cdot a_0=0$ by Theorem \ref{gorenstein condition}.
Thus the map $a_0$ factors through $\delta_{p-q+1}(p)$ by Theorem \ref{gorenstein condition}.
Since $\Hom^{\Z}_R(S(q),S\otimes\bigwedge^{p-q+1}V(p))=0$ by Corollary \ref{negative},
we have $a_0=0$.
By the argument in the proof of (1), we have $a_i=0$ for any $i$.
\end{proof}

Next we shall construct the left isomorphism in \eqref{desired isomorphisms}.

We need the following conventions:
For a positive integer $p$, we put $[p]:=\{1,2,\cdots,p\}$.
When we have an injective map $\sigma:X\to[p]$ from a subset $X$ of $[p]$,
we extend $\sigma$ uniquely to an element $\sigma$ in the symmetric group $\mathfrak{S}_{p}$ by the rule
\begin{itemize}
\item $\sigma(i)<\sigma(j)$ for any $i,j\in[p]\backslash\sigma(X)$ satisfying $i<j$.
\end{itemize}
We have a map
\begin{equation}\label{DV and V}
\bigwedge^qDV\to\Hom_k(\bigwedge^{p+q}V,\bigwedge^pV)
\end{equation}
sending $f=f_q\wedge\cdots\wedge f_2\wedge f_1\in\bigwedge^qDV$ to
\[(v=v_1\wedge v_2\wedge\cdots\wedge v_{p+q}\mapsto
fv=\sum_{\sigma:[q]\to[p+q]}\sgn(\sigma)f_1(v_{\sigma(1)})\cdots f_q(v_{\sigma(q)})v_{\sigma(q+1)}\wedge\cdots\wedge v_{\sigma(p+q)}),\]
where $\sigma$ runs over all injective maps $\sigma:[q]\to[p+q]$, which are uniquely extended to elements in $\mathfrak{S}_{p+q}$.
For the case $p=0$, the above map gives an isomorphism
\begin{equation}\label{DV=DV}
\bigwedge^qDV\simeq D(\bigwedge^qV)
\end{equation}
since $\{y_{i_q}\wedge\cdots\wedge y_{i_1}\ |\ 1\le i_1<\cdots<i_q\le d\}$ is mapped to the dual basis of $\{x_{i_1}\wedge\cdots\wedge x_{i_q}\ |\ 1\le i_1<\cdots<i_q\le d\}$.
The map \eqref{DV and V} gives a map
\begin{eqnarray}\label{from exterior to S}
kG\otimes\bigwedge^qDV&\to&\Hom^{\Z}_R(S\otimes\bigwedge^{p+q}V(p+q),S\otimes\bigwedge^pV(p)),\\
g\otimes f&\mapsto&(s\otimes v\mapsto gs\otimes g(fv)).\nonumber
\end{eqnarray}
For any $0\le q\le d$ the following lemma gives a map
\begin{equation}\label{map from exterior to Koszul}
kG\otimes\bigwedge^qDV\to\Hom_{\C(\modz(R))}(\mathbf{K}(q)[-q],\mathbf{K})
\end{equation}
which is the right map in \eqref{desired isomorphisms}.

\begin{lemma}
For any $0\le q\le d$ and $a\in kG\otimes\bigwedge^qDV$, we have a morphism 
\[\xymatrix@C5em{
&S\otimes\bigwedge^dV(q)\ar[r]^(0.6){(-1)^d\delta_d}\ar[d]^{a}&
\cdots\ar[r]^(0.3){(-1)^{q+2}\delta_{q+2}}&S\otimes\bigwedge^{q+1}V(q)\ar[r]^{(-1)^{q+1}\delta_{q+1}}\ar[d]^{a}&S\otimes\bigwedge^qV(q)\ar[d]^{a}\\
&S\otimes\bigwedge^{d-q}V\ar[r]^(0.6){(-1)^{d-q}\delta_{d-q}}&\cdots\ar[r]^{\delta_2}&S\otimes V\ar[r]^{-\delta_1}&S&&
}\]
of complexes of graded $R$-modules, where the vertical maps are the images of $a$ by \eqref{from exterior to S}.
\end{lemma}

\begin{proof}
We will show that $(-1)^q\delta_pa=a\delta_{p+q}$.
Fix any $a=g\otimes(f_q\wedge\cdots\wedge f_1)\in kG\otimes\bigwedge^qDV$ and $x=s\otimes(v_1\wedge\cdots\wedge v_{p+q})\in S\otimes\bigwedge^{p+q}V(q)$.

Then $\delta_p(a(x))$ is equal to
\begin{equation}\label{E1}
\sum_{\sigma:[q]\to[p+q]}\sum_{\ell=1}^p(-1)^{\ell+1}\sgn(\sigma)g(sv_{\sigma(q+\ell)})\otimes
f_1(v_{\sigma(1)})\cdots f_q(v_{\sigma(q)})(gv_{\sigma(q+1)}\wedge\cdots\check{gv}_{\sigma(q+\ell)}\cdots\wedge gv_{\sigma(p+q)}),
\end{equation}
where $\sigma$ runs over all injective maps $\sigma:[q]\to[p+q]$, which are uniquely extended to elements in $\mathfrak{S}_{p+q}$.
For a pair $(\sigma,\ell)$ appearing in the sum \eqref{E1}, we define an injective map $\tau:[q+1]\to[p+q]$ by putting
\begin{itemize}
\item $\tau(i):=\sigma(i)$ for any $i\in[q]$,
\item $\tau(q+1):=\sigma(q+\ell)$.
\end{itemize}
Since $(-1)^{\ell+1}\sgn(\sigma)=\sgn(\tau)$, the element \eqref{E1} is equal to
\begin{equation}\label{E2}
\sum_{\tau:[q+1]\to[p+q]}\sgn(\tau)g(sv_{\tau(q+1)})\otimes f_1(v_{\tau(1)})\cdots f_q(v_{\tau(q)})(gv_{\tau(q+2)}\wedge\cdots\wedge gv_{\tau(p+q)}),
\end{equation}
where $\tau$ runs over all injective maps $\tau:[q+1]\to[p+q]$, which are uniquely extended to elements in $\mathfrak{S}_{p+q}$.

On the other hand, $a(\delta_{p+q}(x))$ is equal to
\begin{equation}\label{E3}
\sum_{\ell=1}^{p+q}\sum_{\sigma:[q]\to[p+q-1]}(-1)^{\ell+1}\sgn(\sigma)g(sv_\ell)\otimes f_1(v_{\iota_\ell\sigma(1)})\cdots f_q(v_{\iota_\ell\sigma(q)})
(gv_{\iota_\ell\sigma(q+1)}\wedge\cdots\wedge gv_{\iota_\ell\sigma(p+q-1)}),
\end{equation}
where $\sigma$ runs over all injective maps $\sigma:[q]\to[p+q-1]$, which are uniquely extended to elements in $\mathfrak{S}_{p+q-1}$,
and $\iota_\ell$ is the unique bijection $\iota_\ell:[p+q-1]\to[p+q]\backslash\{\ell\}$ which preserves the order on natural numbers.
For a pair $(\ell,\sigma)$ appearing in the sum \eqref{E3}, we define an injective map $\tau:[q+1]\to[p+q]$ by putting
\begin{itemize}
\item $\tau(i):=\iota_\ell\sigma(i)$ for any $i\in[q]$,
\item $\tau(q+1):=\ell$.
\end{itemize}
Since $(-1)^{\ell+1}\sgn(\sigma)=(-1)^q\sgn(\tau)$, the element \eqref{E3} is equal to
\begin{equation}\label{E4}
\sum_{\tau:[q+1]\to[p+q]}(-1)^q\sgn(\tau)g(sv_{\tau(q+1)})\otimes f_1(v_{\tau(1)})\cdots f_q(v_{\tau(q)})(gv_{\tau(q+2)}\wedge\cdots\wedge gv_{\tau(p+q)}),
\end{equation}
where $\tau$ runs over all injective maps $\tau:[q+1]\to[p+q]$, which are uniquely extended to elements in $\mathfrak{S}_{p+q}$.
Comparing \eqref{E2} and \eqref{E4}, we have $(-1)^q\delta_p(a(x))=a(\delta_{p+q}(x))$.
\end{proof}

Using the previous lemma, we have the following result.

\begin{proposition}\label{from exterior to Koszul}
The map \eqref{map from exterior to Koszul} is an isomorphism for any $0\le q\le d$.
\end{proposition}

\begin{proof}
Consider the composition
\[kG\otimes\bigwedge^qDV\xrightarrow{\eqref{map from exterior to Koszul}}\Hom_{\C(\modz(R))}(\mathbf{K}(q)[-q],\mathbf{K})\xrightarrow{\rm rest.}\Hom^{\Z}_R(S\otimes\bigwedge^qV(q),S),\]
where the right map is the restriction to the 0-th terms.
By \eqref{DV=DV} and Corollary \ref{negative}, we have isomorphisms
$\Hom^{\Z}_R(S\otimes\bigwedge^qV(q),S)\simeq\Hom^{\Z}_R(S,S)\otimes\Hom_k(\bigwedge^qV,k)\simeq kG\otimes\bigwedge^qDV$.
Thus the above map is an isomorphism.

Since the map $\Hom_{\C(\modz(R))}(\mathbf{K}(q)[-q],\mathbf{K})\to\Hom^{\Z}_R(\Omega_S^{q+1}k(q),\Omega_Sk)$ is injective by Proposition \ref{Koszul and omega}(1),
the right map is also injective.
Thus we have the assertion.
\end{proof}

%


Now we are ready to prove Theorem \ref{End of tildeU}.

We have desired isomorphisms \eqref{desired isomorphisms} from Propositions \ref{Koszul and omega}(1) and \ref{from exterior to Koszul}.
By Proposition \ref{Koszul and omega}(2), we have isomorphisms
\[\End^{\Z}_R(\widetilde{U})\simeq\bigoplus_{1\le p,q\le d}\Hom^{\Z}_R(\Omega_S^pk(p),\Omega_S^qk(q))\simeq\bigoplus_{1\le p\le q\le d}kG\otimes\bigwedge^{p-q}DV\simeq G*E^{(d)}\]
of $k$-vector spaces.
We have only to check compatibility with the multiplication.

The multiplication in $G*E^{(d)}$ is given by the map
\begin{eqnarray}\label{wedge product}
(kG\otimes\bigwedge^qDV)\otimes(kG\otimes\bigwedge^rDV)&\to&kG\otimes\bigwedge^{q+r}DV,\\
(g\otimes(f_q\wedge\cdots\wedge f_1))\otimes(g'\otimes(f'_r\wedge\cdots\wedge f'_1))&\mapsto&gg'\otimes(f_qg'\wedge\cdots\wedge f_1g'\wedge f'_r\wedge\cdots\wedge f'_1).\nonumber
\end{eqnarray}
We have the following compatibility result.

\begin{lemma}\label{multiplication}
We have a commutative diagram
\[\xymatrix@C=2em{
(kG\otimes\bigwedge^qDV)\otimes\ar[d](kG\otimes\bigwedge^rDV)\ar^{{\rm mult.}}[r]&kG\otimes\bigwedge^{q+r}DV\ar[d]\\
\Hom^{\Z}_R(S\otimes\bigwedge^{p+q}V,S\otimes\bigwedge^pV)\otimes\Hom^{\Z}_R(S\otimes\bigwedge^{p+q+r}V,S\otimes\bigwedge^{p+q}V)
\ar^(0.65){{\rm comp.}}[r]&\Hom^{\Z}_R(S\otimes\bigwedge^{p+q+r}V,S\otimes\bigwedge^pV),
}\]
where the vertical maps are given by \eqref{from exterior to S} and the upper map is given by \eqref{wedge product}.
\end{lemma}

\begin{proof}
Fix $a=g\otimes(f_q\wedge\cdots\wedge f_1)\in kG\otimes\bigwedge^qDV$,
$b=g'\otimes(f'_r\wedge\cdots\wedge f'_1)\in kG\otimes\bigwedge^rDV$ and
$x=s\otimes(v_1\wedge\cdots\wedge v_{p+q+r})\in S\otimes\bigwedge^{p+q+r}V$.

Since $ab=gg'\otimes(f_qg'\wedge\cdots\wedge f_1g'\wedge f'_r\wedge\cdots\wedge f'_1)\in kG\otimes\bigwedge^{q+r}DV$, we have that $(ab)(x)$ equals to
\begin{eqnarray}\label{C1}
\sum_{\mu:[q+r]\to[p+q+r]}&&\sgn(\mu)gg's\otimes f'_1(v_{\mu(1)})\cdots f'_r(v_{\mu(r)})\nonumber\\
&&f_1(g'v_{\mu(1+r)})\cdots f_q(g'v_{\mu(q+r)})
(gg'v_{\mu(1+q+r)}\wedge\cdots\wedge gg'v_{\mu(p+q+r)}),
\end{eqnarray}
where $\mu$ runs over all injective maps $\tau:[q+r]\to[p+q+r]$, which are uniquely extended to elements in $\mathfrak{S}_{q+r}$.

On the other hand, $a(b(x))$ is equal to
\begin{eqnarray}\label{C2}
\sum_{\tau:[r]\to[p+q+r]}&&\sum_{\sigma:[q]\to[p+q]}
\sgn(\tau)\sgn(\sigma)gg's\otimes f'_1(v_{\tau(1)})\cdots f'_r(v_{\tau(r)})\nonumber\\
&&f_1(g'v_{\tau(\sigma(1)+r)})\cdots f_q(g'v_{\tau(\sigma(q)+r)})
(gg'v_{\tau(\sigma(1+q)+r)}\wedge\cdots\wedge gg'v_{\tau(\sigma(p+q)+r)}),
\end{eqnarray}
where $\tau$ runs over all injective maps $\tau:[r]\to[p+q+r]$, which are uniquely extended to elements in $\mathfrak{S}_{p+q+r}$,
and $\sigma$ runs over all injective maps $\sigma:[q]\to[p+q]$, which are uniquely extended to elements in $\mathfrak{S}_{p+q}$.
For a pair $(\tau,\sigma)$ appearing in the sum \eqref{C2}, we define an injective map $\mu:[q+r]\to[p+q+r]$ by putting
\begin{itemize}
\item $\mu(i)=\tau(i)$ for any $i\in[r]$,
\item $\mu(i+r)=\tau(\sigma(i)+r)$ for any $i\in[q]$.
\end{itemize}
Then $\mu$ runs over all injective maps $\mu:[q+r]\to[p+q+r]$.
Since $\sgn(\mu)=\sgn(\tau)\sgn(\sigma)$, we have that \eqref{C2} is equal to \eqref{C1}.
Thus we have $(ab)(x)=a(b(x))$.
\end{proof}

Now we have finished our proof of Theorem \ref{End of tildeU}.
Now Proposition \ref{summand of U} can be shown in a similar way to Proposition \ref{summand of T}.
\qed

\medskip
We shall prove Theorem \ref{End of U}.
Let $e$ and $e'$ be the idempotents of $kG$ defined in \eqref{idempotent}.
The following statements are easily checked from Proposition \ref{summand of U} and Theorem \ref{n-almost split}.
\begin{itemize}
\item[(i)] $[e\widetilde{U}]_{\CM}=R$ and $[e'\widetilde{U}]_{\CM}=e'\widetilde{U}$. 
\item[(ii)] $e'\widetilde{U}$ does not have a non-zero free direct summand.
\end{itemize}
From (i), we have $U=[\widetilde{U}]_{\CM}=R\oplus e'\widetilde{U}$ and $\lend^{\Z}_R(U)\simeq\lend^{\Z}_R(e'\widetilde{U})$. On the other hand,
\[\End^{\Z}_R(e'\widetilde{U})\simeq e'\End^{\Z}_R(\widetilde{U})e'\simeq e'(G*E^{(d)})e'\]
holds by Theorem \ref{End of tildeU}. Thus it remains to show $\End^{\Z}_R(e'\widetilde{U})=\lend^{\Z}_R(e'\widetilde{U})$.
Let $U_p:=\Omega^p_Sk(p)$ for $1\le p\le d$.
Then $e'\widetilde{U}=\bigoplus_{p=1}^de'U_p$, so we have only to show $P^{\Z}(e'U_p,e'U_q)=0$. Assume that there are maps
\[e'U_p\xrightarrow{a}R(r)\xrightarrow{b}e'U_q\]
with $ba\neq0$ for some $r\in\Z$. 
By Theorem \ref{gorenstein condition}, we have a commutative diagram
\[\xymatrix{
e'U_p\ar^a[r]&R(r)\ar^b[r]\ar^{b'}[dr]&e'U_q\\
e'(S\otimes\bigwedge^pV)(p)\ar^{a'}[ur]\ar[u]&&e'(S\otimes\bigwedge^qV)(q)\ar[u]}\]
with $b'a'\neq0$.
Since $\End^{\Z}_R(S)\simeq kG$ by Corollary \ref{negative}, the map $a'$ is a split epimorphism of graded $R$-modules, and so is $a$.
This contradicts to (ii) above. Thus the assertion follows.
\qed

\subsection{Examples: Quivers of the endomorphism algebras}\label{Examples}
Throughout this subsection, we assume that $k$ is an algebraically closed field of characteristic zero.
In this subsection, we give more explicit descriptions of the endomorphism algebras $\underline{\End}^{\Z}_R(T)$ and $\underline{\End}^{\Z}_R(U)$.

Let us start with recalling the quivers of finite dimensional algebras \cite{ARS,ASS}.

\begin{definition}
Let $A$ be a finite dimensional $k$-algebra $A$ and $J_A$ the Jacobson radical of $A$.
Let $\{e_{ij}\ |\ 1\le i\le n,\ 1\le j\le\ell_i\}$ be a complete set of orthogonal primitive idempotents of $A_0$ such that
$A_0e_{ij}\simeq A_0e_{i'j'}$ as $A$-modules if and only if $i=i'$.

The \emph{quiver} $Q$ of $A$ is defined as follows:
The vertices of $Q$ are $1,\cdots,n$.
We draw $d_{ii'}$ arrows from $i$ to $i'$ for $d_{ii'}:=\dim_k(e_{i'j'}(J_A/J_A^2)e_{ij})$, which is independent of $j$ and $j'$.
\end{definition}

\begin{remark}
The above definition is the opposite of \cite{ASS}.
An advantage of our convention is that the directions of morphisms are the same as those of arrows when we consider the quiver of an endomorphism algebra.
\end{remark}

For our cases, the following quivers defined by the representation theory of $G$ are important.

\begin{definition}
Let $G$ be a finite subgroup of $\Gl_d(k)$. Let $V_0=k,V_1,\cdots,V_n$ be the isomorphism classes of simple $kG$-modules. Let $V:=S_1$ be the degree $1$ part of $S$.
\begin{itemize}
\item[(a)] The \emph{McKay quiver} $Q$ of $G$ \cite{A4,Mc,Y} is defined as follows:
The set of vertices is $\{0,1,\cdots,n\}$.
For each vertex $i$, consider the tensor product $V\otimes V_i$ which is a $kG$-module by the diagonal action of $G$. Decompose
\[V\otimes V_i\simeq\bigoplus_{i'=0}^nV_{i'}^{d_{ii'}}\ \ \ (d_{ii'}\in\Z_{\ge0})\]
as an $kG$-module and draw $d_{ii'}$ arrows from $i$ to $i'$.
\item[(b)] We define the \emph{$d$-folded McKay quiver} $Q^{(d)}$ as follows:
The set of vertices is $\{1,2,\cdots,d\}\times\{0,1,\cdots,n\}$.
For any arrow $a:i\to i'$ in $Q$ and $p\in\{1,2,\cdots,d-1\}$, draw an arrow $a:(p,i)\to(p+1,i')$.
\item[(c)] We define the \emph{$d$-folded stable McKay quiver} $\underline{Q}^{(d)}$ by removing the vertices $(p,0)$ for any $1\le p\le d$.
\end{itemize}
\end{definition}

These quivers are especially simple when $G$ is cyclic.

\begin{example}
Let $G$ be a cyclic subgroup of $\Sl_d(k)$ generated by $g={\rm diag}(\zeta^{a_1},\cdots,\zeta^{a_d})$ for a primitive $m$-th root $\zeta$ of unity.
\begin{itemize}
\item[(a)] The isomorphism classes of simple $kG$-modules are $V_i$ ($i\in \Z/m\Z$), where $V_i$ is a one-dimensional $k$-vector space with a basis $\{v_i\}$
such that $gv_i=\zeta^iv_i$.
We have $V_i\otimes V_{i'}\simeq V_{i+i'}$ ($i,i'\in\Z/m\Z$) and $V\simeq V_{a_1}\oplus\cdots\oplus V_{a_d}$ as $kG$-modules.

Thus the McKay quiver $Q$ of $G$ has the vertices $\Z/m\Z$ and the arrows $x_j:i\to i+a_j$ for each $i\in\Z/m\Z$ and $1\le j\le d$.
\item[(b)] The $d$-folded McKay quiver $Q^{(d)}$ of $G$ has the vertices $\Z/m\Z\times\{1,2,\cdots,d\}$ and the arrows $x_j:(p,i)\to(p+1,i+a_j)$ for each $1\le p<d$, $i\in\Z/m\Z$ and $1\le j\le d$.
\item[(c)] The $d$-folded stable McKay quiver $\underline{Q}^{(d)}$ of $G$ is obtained by removing all vertices in $\{1,2,\cdots,d\}\times\{0\}$ from $Q^{(d)}$.
\end{itemize}
\end{example}

It is well-known that the quiver of $\End_R(S)\simeq S*G$ is given by the McKay quiver $Q$ of $G$ \cite{A3,Y}.
Similarly we can draw easily the quivers of $\End^{\Z}_R(T)$, $\End^{\Z}_R(\widetilde{U})$ and $\underline{\End}^{\Z}_R(U)$ as follows:

\begin{proposition}
\begin{itemize}
\item[(a)] The quiver of $\End^{\Z}_R(T)\simeq S^{(d)}*G$ is $Q^{(d)}$.
\item[(b)] The quiver of $\underline{\End}^{\Z}_R(T)\simeq (S^{(d)}*G)/\langle e\rangle$ is $\underline{Q}^{(d)}$.
\item[(c)] The quiver of $\End^{\Z}_R(\widetilde{U})\simeq G*E^{(d)}$ is the opposite quiver $(Q^{(d)})^{\rm op}$ of $Q^{(d)}$.
\end{itemize}
\end{proposition}

\begin{proof}
(a) We have
\begin{eqnarray*}J_{S^{(d)}*G}&=&\left[\begin{array}{ccccccc}
0  &  0&  0&\cdots&0&0&0\\
S_1&  0&  0&\cdots&0&0&0\\
S_2&S_1&  0&\cdots&0&0&0\\
\vdots&\vdots&\vdots&\ddots&\vdots&\vdots\\
S_{d-3}&S_{d-4}&S_{d-5}      &\cdots&0  &0  &  0\\
S_{d-2}&S_{d-3}&S_{d-4}      &\cdots&S_1&0  &  0\\
S_{d-1}&S_{d-2}&S_{d-3}      &\cdots&S_2&S_1&  0
\end{array}\right]\otimes kG,\\
J_{S^{(d)}*G}^2&=&\left[\begin{array}{ccccccc}
0  &  0&  0&\cdots&0&0&0\\
0  &  0&  0&\cdots&0&0&0\\
S_2&  0&  0&\cdots&0&0&0\\
\vdots&\vdots&\vdots&\ddots&\vdots&\vdots\\
S_{d-3}&S_{d-4}&S_{d-5}      &\cdots&0  &0  &  0\\
S_{d-2}&S_{d-3}&S_{d-4}      &\cdots&0  &0  &  0\\
S_{d-1}&S_{d-2}&S_{d-3}      &\cdots&S_2&0  &  0
\end{array}\right]\otimes kG.\end{eqnarray*}
Thus we have
\[S^{(d)}*G/J_{S^{(d)}*G}\simeq k^d\otimes kG\ \mbox{ and }\ J_{S^{(d)}*G}/J_{S^{(d)}*G}^2\simeq V^{d-1}\otimes kG.\]
We have a complete set \eqref{e^p_{ij}} of orthogonal primitive idempotents of $S^{(d)}*G/J_{S^{(d)}*G}$ such that
$(S^{(d)}*G/J_{S^{(d)}*G})e^p_{ij}\simeq (S^{(d)}*G/J_{S^{(d)}*G})e^{p'}_{i'j'}$ if and only if $p=p'$ and $i=i'$.
Thus the vertices of the quiver of $S^{(d)}*G$ is $\{1,2,\cdots,d\}\times\{0,1,\cdots,n\}$.
Let us calculate the number of arrows from $(p,i)$ to $(p',i')$.
Clearly $e^{p'}_{i'j'}(J_{S^{(d)}*G}/J_{S^{(d)}*G}^2)e^{p}_{ij}\neq0$ implies $p'=p+1$.
Moreover, we have
\[e^{p+1}_{i'j'}(J_{S^{(d)}*G}/J_{S^{(d)}*G}^2)e^{p}_{ij}\simeq e_{i'j'}(V\otimes kG)e_{ij}\simeq e_{i'j'}(V\otimes V_{i})\simeq\Hom_{kG}(V_{i'},V\otimes V_{i})\simeq k^{d_{ii'}}.\]
Thus we draw $d_{ii'}$ arrows from $(p,i)$ to $(p+1,i')$, and we have the assertion.

(b) Immediate from (a).

(c) We can prove in a quite similar way to (a).
\end{proof}

The quiver of $\lend^{\Z}_R(U)$ is much more complicated. In the rest of this section, we give some examples.
Let us consider the simplest case $d=2$. Notice that $T\simeq U$ holds in $\lcmz(R)$ in this case.

\begin{example}\label{d=2}
Let $d=2$. Then $\underline{Q}^{(2)}$ is a disjoint union of a Dynkin quiver $\underline{Q}$ and its opposite $\underline{Q}^{\rm op}$.
Moreover $\lend^{\Z}_R(T)\simeq\lend^{\Z}_R(U)$ is isomorphic to the path algebra $k\underline{Q}^{(2)}$.
\end{example}

\begin{proof}
The first assertion is well-known \cite{Mc} (see \cite{A4,Y}).
Since $\underline{Q}^{(2)}$ does not contain a path of length more than one, we have the second assertion.
\end{proof}

Let us consider the case $d=3$. In this case, our assumption that $R$ is an isolated singularity implies that $G$ is cyclic \cite{KN}.
We leave the proof of the following statements to the reader.

\begin{example}
Let $G$ be a cyclic subgroup of $\Sl_3(k)$ generated by ${\rm diag}(\zeta^{a_1},\zeta^{a_2},\zeta^{a_3})$ for a primitive $m$-th root $\zeta$ of unity and $a_1,a_2,a_3\in\Z/m\Z$.
\begin{itemize}
\item[(a)] Assume $a_1=a_2=a_3$. Then the quiver of $\lend^{\Z}_R(U)$ is given in Example \ref{not tilting}.
\item[(b)] Assume $a_1=a_2\neq a_3$. Then the quiver of $\lend^{\Z}_R(U)$ is given by adding an arrow $(a_1,1)\to(-a_1,3)$ to $(\underline{Q}^{(3)})^{\rm op}$ (see Example \ref{pentagon}).
\item[(c)] Assume that $a_1,a_2$ and $a_3$ are mutually distinct. Then the quiver of $\lend^{\Z}_R(U)$ is $(\underline{Q}^{(3)})^{\rm op}$.
\end{itemize}
\end{example}

\begin{example}\label{not tilting}
Let $G$ be a cyclic subgroup of $\Sl_3(k)$ generated by ${\rm diag}(\omega,\omega,\omega)$ for a primitive third root $\omega$ of unity.
The McKay quiver $Q$ of $G$ is the following:
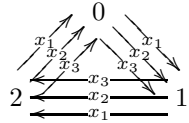
\[\xymatrix@R=2em@C=2em{
&0\ar@<1.5ex>[dr]|{x_1}\ar[dr]|{x_2}\ar@<-1.5ex>[dr]|{x_3}&\\
2\ar@<1.5ex>[ur]|{x_1}\ar[ur]|{x_2}\ar@<-1.5ex>[ur]|{x_3}&&1\ar@<1.5ex>[ll]|{x_1}\ar[ll]|{x_2}\ar@<-1.5ex>[ll]|{x_3}
}\]
The quivers $Q^{(3)}$ and $\underline{Q}^{(3)}$ with relations of $\End_R^{\Z}(T)$ and $\underline{\End}_R^{\Z}(T)$ are the following,
where the vertex $(p,i)$ corresponds to the direct summand $T_{pi}$ defined in \eqref{T_ip}.
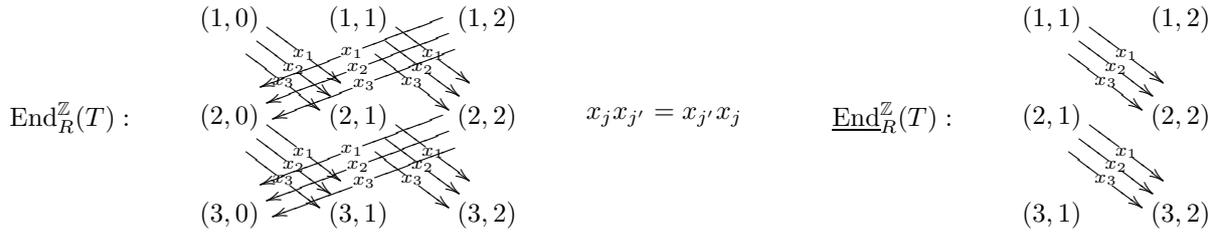
\[\xymatrix@R=2em@C=2em{
&(1,0)\ar@<1.5ex>[dr]|{x_1}\ar[dr]|{x_2}\ar@<-1.5ex>[dr]|{x_3}&(1,1)\ar@<1.5ex>[dr]|{x_1}\ar[dr]|{x_2}\ar@<-1.5ex>[dr]|{x_3}&(1,2)\ar@<1.5ex>[dll]|{x_3}\ar[dll]|{x_2}\ar@<-1.5ex>[dll]|{x_1}\\
\End_R^{\Z}(T):&(2,0)\ar@<1.5ex>[dr]|{x_1}\ar[dr]|{x_2}\ar@<-1.5ex>[dr]|{x_3}&(2,1)\ar@<1.5ex>[dr]|{x_1}\ar[dr]|{x_2}\ar@<-1.5ex>[dr]|{x_3}&(2,2)\ar@<1.5ex>[dll]|{x_3}\ar[dll]|{x_2}\ar@<-1.5ex>[dll]|{x_1}&x_jx_{j'}=x_{j'}x_j\\
&(3,0)&(3,1)&(3,2)
}\ \ \ \ \ \ \ \ \ \ 
\xymatrix@R=2em@C=2em{
&(1,1)\ar@<1.5ex>[dr]|{x_1}\ar[dr]|{x_2}\ar@<-1.5ex>[dr]|{x_3}&(1,2)\\
\underline{\End}_R^{\Z}(T):&(2,1)\ar@<1.5ex>[dr]|{x_1}\ar[dr]|{x_2}\ar@<-1.5ex>[dr]|{x_3}&(2,2)\\
&(3,1)&(3,2)
}\]
On the other hand, the quiver $(Q^{(3)})^{\rm op}$ with relations of $\End_R^{\Z}(\widetilde{U})$ is the following,
where the vertex $(p,i)$ corresponds to the direct summand $U_{pi}$ defined in \eqref{U_ip}.
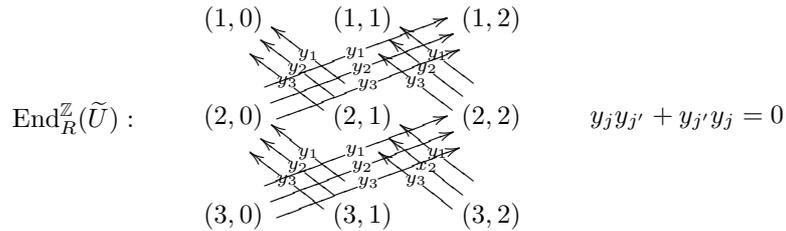
\[\xymatrix@R=2em@C=2em{
&(1,0)&(1,1)&(1,2)\\
\End_R^{\Z}(\widetilde{U}):&(2,0)\ar@<1.5ex>[urr]|{y_1}\ar[urr]|{y_2}\ar@<-1.5ex>[urr]|{y_3}&(2,1)\ar@<1.5ex>[ul]|{y_3}\ar[ul]|{y_2}\ar@<-1.5ex>[ul]|{y_1}
&(2,2)\ar@<1.5ex>[ul]|{y_3}\ar[ul]|{y_2}\ar@<-1.5ex>[ul]|{y_1}&y_jy_{j'}+y_{j'}y_j=0\\
&(3,0)\ar@<1.5ex>[urr]|{y_1}\ar[urr]|{y_2}\ar@<-1.5ex>[urr]|{y_3}&(3,1)\ar@<1.5ex>[ul]|{y_3}\ar[ul]|{y_2}\ar@<-1.5ex>[ul]|{y_1}&(3,2)\ar@<1.5ex>[ul]|{y_3}\ar[ul]|{x_2}\ar@<-1.5ex>[ul]|{y_1}
}\]
The quiver of $\underline{\End}_R^{\Z}(U)$ with relations is the following:
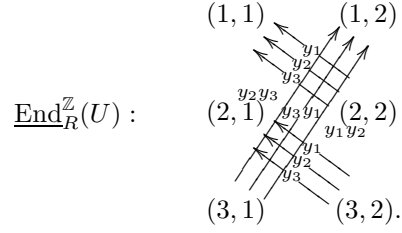
\[\xymatrix@R=2em@C=2em{
&(1,1)&(1,2)\\
\underline{\End}_R^{\Z}(U):&(2,1)&(2,2)\ar@<1.5ex>[ul]|{y_3}\ar[ul]|{y_2}\ar@<-1.5ex>[ul]|{y_1}\\
&(3,1)\ar@<1.5ex>[uur]^{y_2y_3}\ar[uur]|{y_3y_1}\ar@<-1.5ex>[uur]_{y_1y_2}&(3,2)\ar@<1.5ex>[ul]|{y_3}\ar[ul]|{y_2}\ar@<-1.5ex>[ul]|{y_1}.
}\]
In particular for the quiver $\underline{Q}:=\xymatrix{[\bullet\ar@<1ex>[r]\ar[r]\ar@<-1ex>[r]&\bullet]}$, we have a triangle equivalence
\[\lcmz(R)\simeq\D^{\rm b}(\mod k\underline{Q})\times\D^{\rm b}(\mod k\underline{Q})\times\D^{\rm b}(\mod k\underline{Q}).\]
Each factor corresponds to full subcategories $\underline{\CM}^{3\Z+i}(R)$ ($i\in\Z/3\Z$) of $\lcmz(R)$.

Our $R$ gives a non-vanishing example of $(S,S)_i^n$.
Let $S_{\overline{i}}:=\bigoplus_{j\ge0}S_{3j+i}$ for $i=0,1,2$.
For each integer $n$, set
$$
(S_{\overline{i}})_n=
\begin{cases}
S_n & \text{if }n\equiv i\ (\mod 3), \\
0 & \text{otherwise}.
\end{cases}
$$
Then $S_{\overline i}$ is a graded $R$-module.
We have an exact sequence
\[0\to S_{\overline{1}}(-3)\to S_{\overline{2}}(-2)^3\to R(-1)^3\to S_{\overline{1}}\to0\]
by \eqref{Koszul}. Thus we have a triangle
$S_{\overline{1}}(-3)\to S_{\overline{2}}(-2)^3\to S_{\overline{1}}[-1]\to S_{\overline{1}}(-3)[1]$
in $\lcmz(R)$. Applying $\lhom^{\Z}_R(S_{\overline{2}}(-2),-)$ and using $\lhom^{\Z}_R(S,S(-1))=0$, we have
\[\lhom^{\Z}_R(S_{\overline{2}},S_{\overline{1}}(2)[-1])\neq0,\text{ hence }\lhom^{\Z}_R(S,S(2)[-1])\ne 0.\]
In particular, $T=S\oplus S(1)\oplus S(2)$ is not a tilting object.
Note that this does not contradict to Proposition \ref{SS} since $dn+(d-1)i=3\cdot(-1)+2\cdot 2=1>0$.
\end{example}

\begin{example}\label{pentagon}
Let $G$ be a cyclic subgroup of $\Sl_3(k)$ generated by ${\rm diag}(\zeta,\zeta^2,\zeta^2)$ for a primitive fifth root $\zeta$ of unity.
The McKay quiver $Q$ of $G$ is the following:
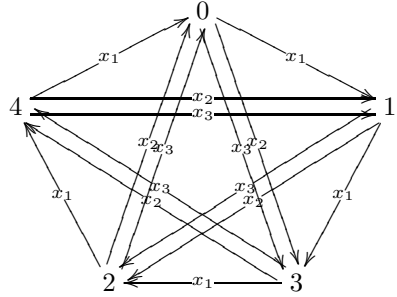
\[\xymatrix{
&&0\ar[rrd]|{x_1}\ar@<.7ex>[rddd]|{x_2}\ar@<-.7ex>[rddd]|{x_3}&&\\
4\ar[urr]|{x_1}\ar@<.7ex>[rrrr]|{x_2}\ar@<-.7ex>[rrrr]|{x_3}&&&&1\ar[ldd]|{x_1}\ar@<.7ex>[llldd]|{x_2}\ar@<-.7ex>[llldd]|{x_3}\\
\\
&2\ar[uul]|{x_1}\ar@<.7ex>[ruuu]|{x_2}\ar@<-.7ex>[ruuu]|{x_3}&&3\ar[ll]|{x_1}\ar@<.7ex>[llluu]|{x_2}\ar@<-.7ex>[llluu]|{x_3}
}\]
The quiver $Q^{(3)}$ of $\End_R^{\Z}(T)$ is the following,
where the vertex $(p,i)$ corresponds to the direct summand $T_{pi}$ defined in \eqref{T_ip}.
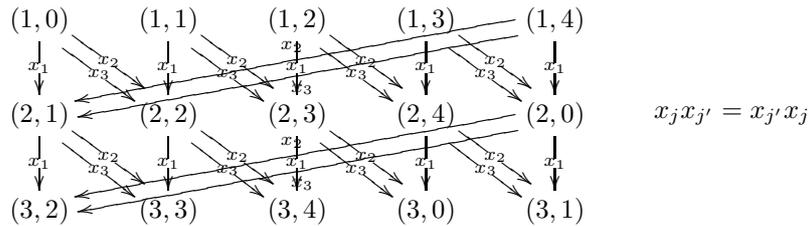
\[\xymatrix@R=2em@C=2em{
(1,0)\ar[d]|{x_1}\ar@<.7ex>[rd]|{x_2}\ar@<-.7ex>[rd]|{x_3}&(1,1)\ar[d]|{x_1}\ar@<.7ex>[rd]|{x_2}\ar@<-.7ex>[rd]|{x_3}&
(1,2)\ar[d]|{x_1}\ar@<.7ex>[rd]|{x_2}\ar@<-.7ex>[rd]|{x_3}&(1,3)\ar[d]|{x_1}\ar@<.7ex>[rd]|{x_2}\ar@<-.7ex>[rd]|{x_3}&(1,4)\ar[d]|{x_1}\ar@<.7ex>[lllld]^{x_3}\ar@<-.7ex>[lllld]_{x_2}\\
(2,1)\ar[d]|{x_1}\ar@<.7ex>[rd]|{x_2}\ar@<-.7ex>[rd]|{x_3}&
(2,2)\ar[d]|{x_1}\ar@<.7ex>[rd]|{x_2}\ar@<-.7ex>[rd]|{x_3}&(2,3)\ar[d]|{x_1}\ar@<.7ex>[rd]|{x_2}\ar@<-.7ex>[rd]|{x_3}&(2,4)\ar[d]|{x_1}\ar@<.7ex>[rd]|{x_2}\ar@<-.7ex>[rd]|{x_3}&
(2,0)\ar[d]|{x_1}\ar@<.7ex>[lllld]^{x_3}\ar@<-.7ex>[lllld]_{x_2}&x_jx_{j'}=x_{j'}x_j\\
(3,2)&(3,3)&(3,4)&(3,0)&(3,1)
}\]
The quiver $\underline{Q}^{(3)}$ of $\underline{\End}_R^{\Z}(T)$ is the following:
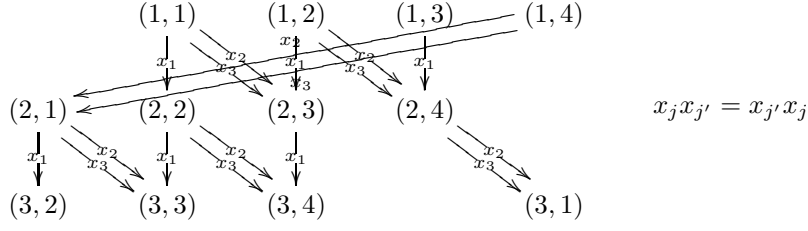
\[\xymatrix@R=2em@C=2em{
&(1,1)\ar[d]|{x_1}\ar@<.7ex>[rd]|{x_2}\ar@<-.7ex>[rd]|{x_3}&
(1,2)\ar[d]|{x_1}\ar@<.7ex>[rd]|{x_2}\ar@<-.7ex>[rd]|{x_3}&(1,3)\ar[d]|{x_1}&(1,4)\ar@<.7ex>[lllld]^{x_3}\ar@<-.7ex>[lllld]_{x_2}\\
(2,1)\ar[d]|{x_1}\ar@<.7ex>[rd]|{x_2}\ar@<-.7ex>[rd]|{x_3}&
(2,2)\ar[d]|{x_1}\ar@<.7ex>[rd]|{x_2}\ar@<-.7ex>[rd]|{x_3}&(2,3)\ar[d]|{x_1}&(2,4)\ar@<.7ex>[rd]|{x_2}\ar@<-.7ex>[rd]|{x_3}&&x_jx_{j'}=x_{j'}x_j\\
(3,2)&(3,3)&(3,4)&&(3,1)
}\]
On the other hand, the quiver $(Q^{(3)})^{\rm op}$ of $\End_R^{\Z}(\widetilde{U})$ is the following,
where the vertex $(p,i)$ corresponds to the direct summand $U_{pi}$ defined in \eqref{U_ip}.
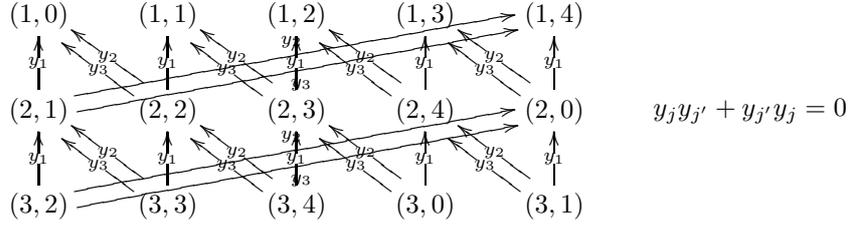
\[\xymatrix@R=2em@C=2em{
(1,0)&(1,1)&(1,2)&(1,3)&(1,4)\\
(2,1)\ar[u]|{y_1}\ar@<.7ex>[rrrru]^{y_2}\ar@<-.7ex>[rrrru]_{y_3}&(2,2)\ar[u]|{y_1}\ar@<.7ex>[ul]|{y_3}\ar@<-.7ex>[ul]|{y_2}&
(2,3)\ar[u]|{y_1}\ar@<.7ex>[ul]|{y_3}\ar@<-.7ex>[ul]|{y_2}&(2,4)\ar[u]|{y_1}\ar@<.7ex>[ul]|{y_3}\ar@<-.7ex>[ul]|{y_2}&(2,0)\ar[u]|{y_1}\ar@<.7ex>[ul]|{y_3}\ar@<-.7ex>[ul]|{y_2}&y_jy_{j'}+y_{j'}y_j=0\\
(3,2)\ar[u]|{y_1}\ar@<.7ex>[rrrru]^{y_2}\ar@<-.7ex>[rrrru]_{y_3}&(3,3)\ar[u]|{y_1}\ar@<.7ex>[ul]|{y_3}\ar@<-.7ex>[ul]|{y_2}&
(3,4)\ar[u]|{y_1}\ar@<.7ex>[ul]|{y_3}\ar@<-.7ex>[ul]|{y_2}&(3,0)\ar[u]|{y_1}\ar@<.7ex>[ul]|{y_3}\ar@<-.7ex>[ul]|{y_2}&(3,1)\ar[u]|{y_1}\ar@<.7ex>[ul]|{y_3}\ar@<-.7ex>[ul]|{y_2}
}\]
The quiver of $\underline{\End}_R^{\Z}(U)$ with relations is the following:
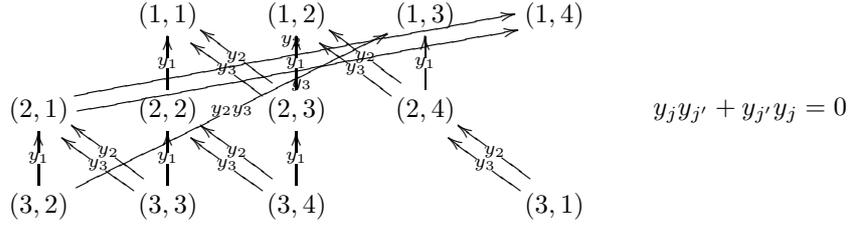
\[\xymatrix@R=2em@C=2em{
&(1,1)&(1,2)&(1,3)&(1,4)\\
(2,1)\ar@<.7ex>[rrrru]^{y_2}\ar@<-.7ex>[rrrru]_{y_3}&(2,2)\ar[u]|{y_1}&
(2,3)\ar[u]|{y_1}\ar@<.7ex>[ul]|{y_3}\ar@<-.7ex>[ul]|{y_2}&(2,4)\ar[u]|{y_1}\ar@<.7ex>[ul]|{y_3}\ar@<-.7ex>[ul]|{y_2}&&y_jy_{j'}+y_{j'}y_j=0\\
(3,2)\ar[u]|{y_1}\ar[rrruu]|{y_2y_3}&(3,3)\ar[u]|{y_1}\ar@<.7ex>[ul]|{y_3}\ar@<-.7ex>[ul]|{y_2}&
(3,4)\ar[u]|{y_1}\ar@<.7ex>[ul]|{y_3}\ar@<-.7ex>[ul]|{y_2}&&(3,1)\ar@<.7ex>[ul]|{y_3}\ar@<-.7ex>[ul]|{y_2}
}\]
\end{example}

\section{Appendix: Algebraic triangulated categories}

In this section we give preliminaries on algebraic triangulated categories.
Let us introduce the following basic notions.

\begin{definition}\cite{Ha,Hel}
Let $\AAA$ be an abelian category and $\BB$ a full subcategory of $\AAA$.
\begin{itemize}
\item[(a)] We say that $\BB$ is \emph{extension-closed} if for any exact sequence $0\to X\to Y\to Z\to 0$ with $X,Z\in\BB$, we have $Y\in\BB$.
In this case, we say that $X\in\BB$ is \emph{relative-projective} if $\Ext^1_{\AAA}(X,\BB)=0$ holds. Similarly we define \emph{relative-injective} objects in $\BB$.
\item[(b)] An extension closed subcategory $\BB$ of an abelian category $\AAA$ is called \emph{Frobenius} if the following conditions are satisfied:
\begin{itemize}
\item[(i)] An object in $\BB$ is relative-projective if and only if it is relative-injective.
\item[(ii)] For any $X\in\BB$, there exist exact sequences $0\to Y\to P\to X\to0$ and $0\to X\to I\to Z\to0$ in $\AAA$ such that
$P\in\BB$ is relative-projective, $I\in\BB$ is relative-injective and $Y,Z\in\BB$.
\end{itemize}
\item[(c)] For a Frobenius category $\BB$, we define the \emph{stable category} $\underline{\BB}$ as follows:
The objects of $\underline{\BB}$ are the same as $\BB$, and the morphism set is given by
\[\underline{\Hom}_{\BB}(X,Y):=\Hom_{\BB}(X,Y)/P(X,Y)\]
for any $X,Y\in\BB$, where $P(X,Y)$ is the submodule of $\Hom_{\BB}(X,Y)$ consisting of morphisms which factor through relative-projective objects in $\BB$.
\end{itemize}
\end{definition}

We refer to \cite{Ke2} for an axiomatic definition of a Frobenius category, which is slightly more general when $\BB$ is not small.

An important property of Frobenius categories found by Happel is the following.

\begin{definition-theorem}\cite{Ha}\cite{Ke2}
The stable category $\underline{\BB}$ of a Frobenius category $\BB$ has a structure of a triangulated category.
Such a triangulated category is called \emph{algebraic}.
\end{definition-theorem}

One of the advantages of algebraic triangulated categories is that we can realize them as homotopy categories:
Let $\BB$ be a Frobenius category, and $\PP$ the full subcategory of relative-projective objects in $\BB$.
Let $\C^{\rm ac}(\PP)$ be the category of chain complexes over $\PP$ which are obtained by gluing short exact sequences in $\BB$.
As usual, the homotopy category $\K^{\rm ac}(\PP)$ has a structure of a triangulated category.

\begin{proposition}\label{from stable to homotopy}
We have a triangle equivalence
\[Z^0:\K^{\rm ac}(\PP)\to\underline{\BB}.\]
\end{proposition}

\subsection{Proof of Theorem \ref{equivalence}}

The tilting theorem for algebraic triangulated categories was given by Keller \cite[(4.3)]{Ke1} (see also \cite[(6.5)]{Kr}),
and its weakest form is Theorem \ref{equivalence}.
We include a proof of Theorem \ref{equivalence} for the convenience of the reader.

First we need the following well-known observation.

\begin{lemma}\label{devissage}
Let $F:\TT\to\TT^\prime$ be a triangle functor of triangulated categories and $U\in\TT$ an object.
If $F_{U,U[n]}:\Hom_{\TT}(U,U[n])\to\Hom_{\TT^\prime}(FU, FU[n])$ is an isomorphism for any $n\in\Z$, then $F:\mathsf{thick}_{\TT}(U)\to\TT'$ is fully faithful.
\end{lemma}

Next we need a general observation on chain complexes.

For an additive category $\PP$, we denote by $\C(\PP)$ (respectively, $\K(\PP)$) the category (respectively, homotopy category) of chain complexes over $\PP$.
For two complexes $X=(\cdots\to X^n\stackrel{d_X^n}{\to}X^{n+1}\to\cdots)$ and $Y=(\cdots\to Y^n\stackrel{d_Y^n}{\to}Y^{n+1}\to\cdots)$ in $\C(\PP)$,
we have a complex
\[\HHom(X,Y)=(\cdots\to \HHom(X,Y)^n\stackrel{d^n}{\to}\HHom(X,Y)^{n+1}\to\cdots)\]
where
\[\HHom(X,Y)^n:=\prod_{p\in\Z}\Hom_{\PP}(X^p,Y^{p+n})\]
and the differential is given by
\[d^n((\phi^p)_{p\in\Z})=(d_Y^{p+n}\circ\phi^p-(-1)^n\phi^{p+1}\circ d_X^p)_{p\in\Z}.\]
In particular $\HHom(X,X)$ has a structure of a differential graded (=DG) algebra.
It is easy to check
\[\Hom_{\C(\PP)}(X,Y)\simeq Z^0(\HHom(X,Y)))\ \mbox{ and }\ \Hom_{\K(\PP)}(X,Y)\simeq H^0(\HHom(X,Y))).\]

For a DG algebra $A$, we denote by $\C A$ (respectively, $\K A$, $\D A$)
the category (respectively, homotopy category, derived category) of DG $A$-modules (see \cite{Ke1,Kr}).
As usual we denote the subcategory $\mathsf{thick}_{\D A}(A)$ by $\per A$.
If $A$ is concentrated in degree $0$, then $\per A$ coincides with $\K^{\bo}(\proj A)$.

\begin{lemma}\label{homotopy to DG}
For any complex $U=(\cdots\to U^n\stackrel{d_U^n}{\to}U^{n+1}\to\cdots)\in\C(\PP)$, define a DG algebra by $A:=\HHom(U,U)$.
Then we have a triangle equivalence $\mathsf{thick}_{\K(\PP)}(U)\to\per A$ up to direct summands.
\end{lemma}

\begin{proof}
We have a functor
\[\C(\PP)\to\C A,\ X\mapsto\HHom(U,X).\]
Since this functor sends a null-homotopic morphism of complexes over $\PP$ to a null-homotopic morphism of DG $A$-modules,
we have a triangle functor $\K(\PP)\to\K A$.
Composing with the canonical functor $\K A\to \D A$, we have a triangle functor
$\K(\PP)\to\D A$.
Since $U$ is sent to $A$, this induces a triangle functor
\[F:\mathsf{thick}_{\K(\PP)}(U)\to\mathsf{thick}_{\D A}(A)=\per A.\]
Since we have a commutative diagram
\[\xymatrix@C=4em{\Hom_{\K(\PP)}(U,U[n])\ar^{F_{U,U[n]}}[r]\ar[d]^{\wr}&\Hom_{\D A}(A,A[n])\ar[d]^{\wr}\\
H^n(A)\ar@{=}[r]&H^n(A)
}\]
for any $n\in\Z$, the functor $F$ is a triangle equivalence up to direct summands by Lemma \ref{devissage}.
Thus we have the desired equivalence.
\end{proof}

Finally we need the following observation on quasi-isomorphisms of DG algebras.

\begin{lemma}\label{quasi-isomorphism}
Let $f:B\to A$ be a quasi-isomorphism of DG algebras. Then we have a triangle equivalence $\per A\to\per B$. 
\end{lemma}

\begin{proof}
Although this is elementary, we give a proof for the convenience of the reader.

Since any DG $A$-module (respectively, morphism of DG $A$-modules) can be regarded as a DG $B$-module (respectively, morphism of DG $B$-modules),
we have a functor $\C A\to\C B$.
Since any null-homotopic morphism of DG $A$-modules is also null-homotopic as a morphism of DG $B$-modules,
we have an induced functor $\K A\to \K B$.
Since any quasi-isomorphism of DG $A$-modules is also a quasi-isomorphism of DG $B$-modules,
we get an induced functor $G:\D A\to \D B$.

Since $f:B\to A$ is an isomorphism in $\D B$, we have an isomorphism $f[n]^{-1}\cdot f:\Hom_{\D B}(A,A[n])\to\Hom_{\D B}(B,B[n])$ for any $n\in\Z$.
Since we have a commutative diagram
\[\xymatrix@C=4em{
\Hom_{\D A}(A,A[n])\ar[d]^{\wr}\ar[r]^{G_{A,A[n]}}&\Hom_{\D B}(A,A[n])\ar[r]_{\sim}^{f[n]^{-1}\cdot f}&\Hom_{\D B}(B,B[n])\ar[d]^{\wr}\\
H^n(A)&&H^n(B)\ar[ll]^{\sim}_{H^n(f)}}\]
the map $G_{A,A[n]}:\Hom_{\D A}(A,A[n])\to\Hom_{\D B}(A,A[n])$ is an isomorphism for any $n\in\Z$.
Thus the functor $\per A\to\D B$ is fully faithful by Lemma \ref{devissage}.
Since $A\simeq B$ in $\D B$ by our assumption, we obtain a triangle equivalence $\per A\to \per B$ up to direct summands.
This is dense since we have $X\simeq A\otimes_B^{\bf L}X$ in $\D A$ for any $X\in\D B$.
\end{proof}

Now we are ready to prove Theorem \ref{equivalence}.

We assume that $\TT=\underline{\BB}$ for a Frobenius category $\BB$.
Without loss of generality, we can assume $\TT=\K^{\rm ac}(\PP)$ by Proposition \ref{from stable to homotopy}.
For the tilting object $U\in\TT=\K^{\rm ac}(\PP)$, define a DG algebra by $A:=\HHom(U,U)$.
By Lemma \ref{homotopy to DG}, we have a triangle equivalence
\begin{equation}\label{equivalence1}
\TT=\mathsf{thick}_{\TT}(U)\xrightarrow{}\per A
\end{equation}
up to direct summands. Thus we have
\begin{equation}\label{homology of A}
H^n(A)\simeq\Hom_{\D A}(A,A[n])\simeq\Hom_{\K(\PP)}(U,U[n])=
\left\{\begin{array}{cc}
0&n\neq0,\\
\End_{\TT}(U)&n=0.
\end{array}\right.
\end{equation}

Now we denote by $B$ the DG subalgebra of $A$ defined by
\[B^n:=\left\{\begin{array}{cc}
A^n&n<0,\\
\Ker d_A^0&n=0,\\
0&n>0.
\end{array}\right.\]
By \eqref{homology of A}, the natural inclusion $B\to A$ is a quasi-isomorphism of DG algebras,
and the natural surjection $B^0=\Ker d_A^0\to H^0(A)$ induces a quasi-isomorphism $B\to H^0(A)$ of DG algebras,
where we regard $H^0(A)$ as a DG algebra concentrated in degree $0$.
By Lemma \ref{quasi-isomorphism}, we have triangle equivalences
\begin{equation}\label{equivalence2}
\per A\xrightarrow{} \per B\xleftarrow{}\per H^0(A).
\end{equation}
Composing \eqref{equivalence1} and \eqref{equivalence2}, we have the desired triangle equivalence
\[\TT\to\per H^0(A)=\K^{\bo}(\proj\End_{\TT}(U))\]
up to direct summands.
\qed


\begin{thebibliography}{13}
\bibitem[AI]{AI} T. Aihara, O. Iyama, \emph{Silting mutation in triangulated categories}, arXiv:1009.3370.
\bibitem[AIR]{AIR} C. Amiot, O. Iyama, I. Reiten, \emph{Stable categories of Cohen-Macaulay modules and cluster categories}, in preparation.
\bibitem[Ar]{Ar} T. Araya, \emph{Exceptional sequences over graded Cohen-Macaulay rings}, Math. J. Okayama Univ.  41  (1999), 81--102 (2001).
\bibitem[A1]{A1} M. Auslander, \emph{On the purity of the branch locus}, Amer. J. Math. 84 1962 116--125. 
\bibitem[A2]{A2} M. Auslander, \emph{Coherent functors},
1966 Proc. Conf. Categorical Algebra (La Jolla, Calif., 1965) pp. 189--231 Springer, New York.
\bibitem[A3]{A3} M. Auslander, \emph{Functors and morphisms determined by objects}, Representation theory of algebras (Proc. Conf., Temple Univ., Philadelphia, Pa., 1976), pp. 1--244. Lecture Notes in Pure Appl. Math., Vol. 37, Dekker, New York, 1978. 
\bibitem[A4]{A4} M. Auslander, \emph{Rational singularities and almost split sequences},
Trans. Amer. Math. Soc. 293 (1986), no. 2, 511--531. 
\bibitem[ABr]{ABr} M. Auslander, M. Bridger, \emph{Stable module theory}, Memoirs of the American Mathematical Society, No. 94 American Mathematical Society, Providence, R.I. 1969.
\bibitem[AG]{AG} M. Auslander, O. Goldman, \emph{The Brauer group of a commutative ring},  Trans. Amer. Math. Soc.  97  1960 367--409.
\bibitem[AR1]{AR1} M. Auslander, I. Reiten, \emph{Almost split sequences for $Z$-graded rings},
Singularities, representation of algebras, and vector bundles (Lambrecht, 1985),  232--243, Lecture Notes in Math., 1273, Springer, Berlin, 1987. 
\bibitem[AR2]{AR2} M. Auslander, I. Reiten, \emph{Cohen-Macaulay modules for graded Cohen-Macaulay rings and their completions},
Commutative algebra (Berkeley, CA, 1987), 21--31, Math. Sci. Res. Inst. Publ., 15, Springer, New York, 1989. 
\bibitem[ARS]{ARS} M. Auslander, I. Reiten, S. O. Smal\o, \emph{Representation theory of Artin algebras}, Cambridge Studies in Advanced Mathematics, 36.
Cambridge University Press, Cambridge, 1997.
\bibitem[ASS]{ASS} I. Assem, D. Simson, A. Skowro\'{n}ski, \emph{Elements of the representation theory of associative algebras, Vol. 1}, Techniques of representation theory. London Mathematical Society Student Texts, 65.
Cambridge University Press, Cambridge, 2006.
\bibitem[BR]{BR} H. Bass, A. Roy, \emph{Lectures on topics in algebraic $K$-theory},
Tata Institute of Fundamental Research Lectures on Mathematics, No. 41 Tata Institute of Fundamental Research, Bombay 1967.
\bibitem[Bo]{Bo} A. I. Bondal, \emph{Representations of associative algebras and coherent sheaves},
Izv. Akad. Nauk SSSR Ser. Mat. 53 (1989), no. 1, 25--44; translation in Math. USSR-Izv. 34 (1990), no. 1, 23--42
\bibitem[BK]{BK} A. I. Bondal, M. M. Kapranov, \emph{Representable functors, Serre functors, and reconstructions},
Izv. Akad. Nauk SSSR Ser. Mat. 53 (1989), no. 6, 1183--1205, 1337; translation in Math. USSR-Izv. 35 (1990), no. 3, 519--541.
\bibitem[BH]{BH} W. Bruns, J. Herzog, \emph{Cohen-Macaulay rings (revised edition)}, Cambridge Studies in Advanced Mathematics, 39, Cambridge University Press, Cambridge, 1998.
\bibitem[Bu]{Bu} R.-O. Buchweitz, \emph{Morita contexts, idempotents, and Hochschild cohomology---with applications to invariant rings},
Commutative algebra (Grenoble/Lyon, 2001), 25--53, Contemp. Math., 331, Amer. Math. Soc., Providence, RI, 2003.
\bibitem[BIKR]{BIKR} I. Burban, O. Iyama, B. Keller, I. Reiten, \emph{Cluster tilting for one-dimensional hypersurface singularities}, Adv. Math. 217 (2008), no. 6, 2443--2484.
\bibitem[CR1]{CR1} C. W. Curtis, I. Reiner, \emph{Methods of representation theory. Vol. I. With applications to finite groups and orders},
Pure and Applied Mathematics. A Wiley-Interscience Publication. John Wiley \& Sons, Inc., New York, 1981.
\bibitem[CR2]{CR2} C. W. Curtis, I. Reiner, \emph{Methods of representation theory. Vol. II. With applications to finite groups and orders},
Pure and Applied Mathematics (New York). A Wiley-Interscience Publication. John Wiley \& Sons, Inc., New York, 1987.
\bibitem[DH]{DH} H. Dao, C. Huneke, \emph{Vanishing of Ext, cluster tilting modules and finite global dimension of endomorphism rings}, arXiv:1005.5359.
\bibitem[Ha]{Ha} D. Happel, \emph{Triangulated categories in the representation theory of finite-dimensional algebras},
London Mathematical Society Lecture Note Series, 119. Cambridge University Press, Cambridge, 1988.
\bibitem[Hel]{Hel} A. Heller, \emph{The loop-space functor in homological algebra},
Trans. Amer. Math. Soc. 96 (1960), 382--394.
\bibitem[Her]{Her} J. Herzog, \emph{Ringe mit nur endlich vielen Isomorphieklassen von maximalen, unzerlegbaren Cohen-Macaulay-Moduln}, Math. Ann. 233 (1978), no. 1, 21--34.
\bibitem[I1]{I1} O. Iyama, \emph{Higher-dimensional Auslander-Reiten theory on maximal orthogonal subcategories}, Adv. Math. 210 (2007), no. 1, 22--50.
\bibitem[I2]{I2} O. Iyama, \emph{Auslander correspondence}, Adv. Math. 210 (2007), no. 1, 51--82.
\bibitem[IR]{IR} O. Iyama, I. Reiten, \emph{Fomin-Zelevinsky mutation and tilting modules over Calabi-Yau algebras}, Amer. J. Math. 130 (2008), no. 4, 1087--1149.
\bibitem[IW]{IW} O. Iyama, M. Wemyss, \emph{Auslander-Reiten Duality and Maximal Modifications for Non-isolated Singularities}, arXiv:1007.1296.
\bibitem[IY]{IY} O. Iyama, Y. Yoshino, \emph{Mutation in triangulated categories and rigid Cohen-Macaulay modules}, Invent. Math. 172 (2008), no. 1, 117--168.
\bibitem[KST1]{KST1} H. Kajiura, K. Saito, A. Takahashi, \emph{Matrix factorization and representations of quivers. II. Type $ADE$ case}, Adv. Math. 211 (2007), no. 1, 327--362.
\bibitem[KST2]{KST2} H. Kajiura, K. Saito, A. Takahashi, \emph{Triangulated categories of matrix factorizations for regular systems of weights with $\epsilon=-1$}, Adv. Math. 220 (2009), no. 5, 1602--1654.
\bibitem[Ke1]{Ke1} B. Keller, \emph{Deriving DG categories}, Ann. Sci. Ecole Norm. Sup. (4)  27  (1994),  no. 1, 63--102.
\bibitem[Ke2]{Ke2} B. Keller, \emph{On differential graded categories}, International Congress of Mathematicians. Vol. II,  151--190, Eur. Math. Soc., Zuich, 2006.
\bibitem[KR1]{KR1} B. Keller, I. Reiten, \emph{Cluster-tilted algebras are Gorenstein and stably Calabi-Yau},  Adv. Math.  211  (2007),  no. 1, 123--151.
\bibitem[KR2]{KR2} B. Keller, I. Reiten, \emph{Acyclic Calabi-Yau categories. With an appendix by Michel Van den Bergh},  Compos. Math.  144  (2008),  no. 5, 1332--1348.
\bibitem[KMV]{KMV} B. Keller, D. Murfet, M. Van den Bergh, \emph{On two examples by Iyama and Yoshino}, to appear in Compos. Math., arXiv:0803.0720.
\bibitem[KV]{KV} B. Keller, D. Vossieck, \emph{Aisles in derived categories},
Deuxieme Contact Franco-Belge en Algebre (Faulx-les-Tombes, 1987).
Bull. Soc. Math. Belg. Ser. A 40 (1988), no. 2, 239--253.
\bibitem[Kr]{Kr} H. Krause, \emph{Derived categories, resolutions, and Brown representability and Exercises}, presentation at Summer School Chicago 2004, available at http://www2.math.uni-paderborn.de/people/henning-krause.html
\bibitem[KN]{KN} K. Kurano, S. Nishi, \emph{Gorenstein isolated quotient singularities of odd prime dimension are cyclic}, arXiv:0903.3270.
\bibitem[LP]{LP} H. Lenzing, J. A. de la Pena, \emph{Extended canonical algebras and Fuchsian singularities}, to appear in Math. Z., arXiv:math/0611532.
\bibitem[Ma]{Ma} H. Matsumura, \emph{Commutative ring theory}, Translated from the Japanese by M. Reid, Second edition, Cambridge Studies in Advanced Mathematics, 8, 
Cambridge University Press, Cambridge, 1989.
\bibitem[Mc]{Mc} J. McKay, \emph{Graphs, singularities, and finite groups},  The Santa Cruz Conference on Finite Groups (Univ. California, Santa Cruz, Calif., 1979),
pp. 183--186, Proc. Sympos. Pure Math., 37, Amer. Math. Soc., Providence, R.I., 1980.
\bibitem[O]{O} D. Orlov, \emph{Derived categories of coherent sheaves and triangulated categories of singularities}, Algebra, arithmetic, and geometry: in honor of Yu. I. Manin. Vol. II, 503--531, Progr. Math., 270, Birkhauser Boston, Inc., Boston, MA, 2009.
\bibitem[R]{R} J. Rickard, \emph{Morita theory for derived categories}, J. London Math. Soc. (2) 39 (1989), no. 3, 436-456.
\bibitem[S]{S} J.-P. Serre, \emph{Local fields}, Graduate Texts in Mathematics, 67. Springer-Verlag, New York-Berlin, 1979.
\bibitem[T]{T} A. Takahashi, \emph{Matrix Factorizations and Representations of Quivers I}, arXiv:math/0506347.
\bibitem[TV]{TV} L. de Thanhoffer de Volcsey, M. Van den Bergh, \emph{Explicit models for some stable categories of maximal Cohen-Macaulay modules}, arXiv:1006.2021.
\bibitem[U]{U} K. Ueda, \emph{Triangulated categories of Gorenstein cyclic quotient singularities}, Proc. Amer. Math. Soc. 136 (2008), no. 8, 2745--2747.
\bibitem[Y]{Y} Y. Yoshino, \emph{Cohen-Macaulay modules over Cohen-Macaulay rings}, London Mathematical Society Lecture Note Series, 146, Cambridge University Press, Cambridge, 1990.
\end{thebibliography}
\end{document}